\documentclass[12pt]{amsart}
\usepackage[margin=1in]{geometry}
\usepackage{amsmath,latexsym,amssymb,amsthm,graphicx}
\usepackage{mathtools}
\usepackage{subfig}
\usepackage{xfrac}
\usepackage{fix-cm}
\usepackage{color}
\usepackage{verbatim}
\DeclareMathOperator\arctanh{arctanh}
\usepackage[bookmarks=true,%
    colorlinks=true,%
    linkcolor=blue,%
    citecolor=blue,%
    filecolor=blue,%
    menucolor=blue,%
    urlcolor=blue]{hyperref}
\usepackage{doi}
\usepackage{stmaryrd}
\usepackage{multirow}
\usepackage{float}
\usepackage{tabularx}
\usepackage[shortlabels]{enumitem}
\usepackage{pdfsync}
\usepackage{tikz}
\usepackage{tikz-cd}
\usepackage{float}
\usepackage[utf8]{inputenc}
\usepackage[textsize=small]{todonotes}

\usepackage{hyphenat}
\usepackage{booktabs}

\graphicspath{{mpdraws/}{draws/}}

\setenumerate[1]{label=(\arabic*)}

\numberwithin{equation}{section}
\numberwithin{figure}{section}
\newtheorem{mainthm}{Theorem}
\newtheorem*{theorem*}{Theorem}

\newtheorem{theorem}{Theorem}[section]
\newtheorem{lemma}[theorem]{Lemma}

\newtheorem{corollary}[theorem]{Corollary}
\newtheorem{proposition}[theorem]{Proposition}

\newtheorem*{ex*}{Exercise}

\theoremstyle{definition}
\newtheorem{remark}[theorem]{Remark}

\newtheorem{definition}[theorem]{Definition}
\newtheorem{question}[theorem]{Question}

\newtheorem{conjecture}[theorem]{Conjecture}

\newtheorem*{acknowledgements}{Acknowledgements}

\DeclareMathOperator{\Teich}{Teich}

\DeclareMathOperator{\len}{len}

\DeclareMathOperator{\EL}{EL}
\DeclareMathOperator{\sys}{sys}
\DeclareMathOperator{\Area}{Area}

\renewcommand{\epsilon}{\varepsilon}

\DeclareMathOperator{\PSL}{\mathit{PSL}}

\DeclareMathOperator{\arcsinh}{arcsinh}
\DeclareMathOperator{\arccosh}{arccosh}

\newcommand{\VR}{{\rm V}_{\rm R}}

\usetikzlibrary{arrows}
\tikzset{
    labl/.style={anchor=south, rotate=90, inner sep=.5mm}
}
\tikzstyle{every picture}=[> = to]
\tikzset{cdlabel/.style={execute at begin node=$\scriptstyle,execute at end node=$}}
\tikzset{implication/.style={double equal sign distance, -implies}}
\tikzset{biimplication/.style={double equal sign distance, implies-implies}}

\begin{document}
\title{Comparing hyperbolic and extremal lengths for shortest curves}
\author[Martínez-Granado]{Dídac Martínez-Granado}
\address{Department of Mathematics\\University of Luxembourg\\Av. de la Fonte 6, Esch-sur-Alzette, L-4364, Luxembourg}       
\email{didac.martinezgranado@uni.lu}

\author{Franco Vargas Pallete}
\address{IMPA - Instituto Nacional de Matem\'atica Pura e Aplicada, Rio de Janeiro, RJ, Brasil, 22460-320}
\email{vargas.pallete@impa.br}

\begin{abstract}
We give a lower bound for the widths of the collars of certain short partial pants 
decomposition of the surface. Then we apply this to obtain upper bounds of the renormalized volume of certain Schottky manifolds in terms of the hyperbolic length of compressible curves.
\end{abstract}
\setcounter{tocdepth}{1}
\maketitle

\section{Introduction}

The geometry of surfaces admits two intertwined viewpoints: the hyperbolic
geometry arising from the uniformization theorem, and the conformal geometry
encoded by extremal length. Each provides a natural notion of curve length.
For a hyperbolic surface $\Sigma$, the \emph{hyperbolic length} of a curve is
the length of its geodesic representative, while the \emph{extremal length} is
a conformal invariant defined via a variational principle
(see Definition~\ref{def:EL} for a precise definition). Extremal length plays a
central role in Teichmüller theory
(\cite{Ker90:Asymptotics, LR11:Quasiconvex, Min96:ExtremalLength}), but unlike
hyperbolic length it is often difficult to compute directly.

A recurring theme is to compare these two notions of length. For instance,
Maskit’s comparison theorem \cite{Maskit85:ComparisonLengths}, refined by Buser and
Sarnak \cite{BS94:PeriodMatrix}, shows that the extremal length of the
\emph{hyperbolic systole}—the shortest nontrivial closed geodesic—is bounded
above by a universal multiple of its hyperbolic length. In particular, this
controls the \emph{extremal sytole} by the hyperbolic systole.

In this paper we extend such comparisons beyond systoles, to curves appearing
in \emph{short partial pants decompositions}. These are obtained by
successively cutting along the hyperbolic systole and then selecting the next
shortest closed geodesic, up to $n < 3g-3$ disjoint curves in total. We further assume
that each geodesic has hyperbolic length at most $\epsilon$, for a universal
constant $\epsilon > 0$ defined below. We shall refer to these as \emph{ $\epsilon$-short partial pants decompositions.}

Our result shows that even in this broader setting, the complementary
components of the decomposition still contain curves that enjoy collars of
uniform width. This result is an extension of Buser and
Sarnak's result for systoles, stated as Theorem~\ref{thm:busersarnak}. This geometric control yields a comparison between
extremal and hyperbolic length:

\begin{mainthm}[Main Theorem, informal]
On any closed hyperbolic surface, 
for every $\epsilon$-short partial pants decomposition, the shortest curve in each component of the complement admits a uniform width (in terms of $\epsilon$) embedded collar. Consequently, its extremal length is bounded above by
a universal multiple of its hyperbolic length.
\end{mainthm}

The precise quantitative version, with explicit constants, appears in
Theorem~\ref{thm:mainthm}. The inequality between extremal and hyperbolic
length then follows directly from a classical comparison due to Maskit (see
Theorem~\ref{thm:upperEL} in the Appendix). The main novelty of our approach is a
uniform lower bound on the width of half collars around the next shortest curve
in an $\epsilon$-short partial pants decomposition. Our argument is inspired by
Buser--Sarnak \cite[Section~4]{BS94:PeriodMatrix}, but goes beyond their
methods, since their results alone do not yield the required uniform width
estimate. Indeed, given any $\epsilon$-short partial pants decomposition, our result applies to the next shortest curve on each component of their complement regardless of their topological type. Moreover, such examples exist for any topological type of a partial pants decomposition on a surface, by making the lengths of the curves of the desired topological type small enough in Fenchel-Nielsen coordinates.

As an application, we obtain upper bounds for the \emph{renormalized volume} of
certain Schottky manifolds in terms of the hyperbolic length. The renormalized
volume $\VR$~\cite{KS08} is a geometric invariant that assigns a finite notion of volume to
hyperbolic 3--manifolds that otherwise have infinite hyperbolic volume (see
Section~\ref{subsec:renormalizedvolume} for precise
definitions). In the statement below, we say a collection of disjoint simple
curves is \emph{suitable} if each component of the complement has genus
$1$ or $2$, and any component of genus $2$ has exactly one boundary component. Recall that an essential simple curve is in the boundary of a 3-manifold $M^3$ is  \emph{compressible} if it bounds a disk in $M^3$, and $A_w(\gamma)$ denotes a \emph{half collar} of width $w$ about the essential simple curve $\gamma$ (see Section~\ref{subsec:systoles} for a definition).
\begin{mainthm}
Let $w_0,\epsilon,c$ be the constants of Theorem~\ref{thm:mainthm}. Let $\Sigma$ be
a closed Riemann surface and let $\mathcal{C} = \lbrace \gamma_k\rbrace_{1\leq
k \leq n}$ be an $\epsilon$-short partial pants decomposition, each admitting an embedded half collar $A_{w_0}(\gamma_k)$ in $\Sigma$,
and suppose that $\mathcal{C}$ is suitable. Then we can extend $\mathcal{C}$ to
a collection $\mathcal{C}'$ of disjoint curves so that each component of
$\Sigma\setminus\mathcal{C}'$ has genus $1$. For any Schottky manifold $M^3$
where all curves in $\mathcal{C}'$ are compressible, we obtain the bound
\[
    \VR(M) \leq \pi(g-1) +
    c\left(\sum_{i=1}^{g-1} \sqrt{\ell(\gamma_{k_i}; \Sigma)} \right)^2.
\]

Moreover, if
\[
    \sqrt{c} \left(\sum_{i=1}^{g-1} \sqrt{\ell(\gamma_{k_i}; \Sigma)} \right) \leq
    \sqrt{\pi(g-1)},
\]
then
\[
    \VR(M) \leq \pi(g-1)\left(3-\frac{\pi(g-1)}{c\left(\sum_{i=1}^{g-1}
    \sqrt{\ell(\gamma_{k_i}; \Sigma)} \right)^2}\right).
\]
In particular, $\VR \leq 0$ whenever
\[
    \sqrt{c}\sum_{i=1}^{g-1} \sqrt{\ell(\gamma_{k_i}; \Sigma)} \leq
    \sqrt{\tfrac{\pi(g-1)}{3}}.
\]
\label{mainthm:renormalized}
\end{mainthm}

Upper bounds for the renormalized volume in terms of the hyperbolic length are useful to understand the growth of $\VR$ as a function of the genus. Moreover, hyperbolic length of a curve is typically easier to estimate than extremal length, which is difficult to estimate explicitly in general. A similar result appears in work of Cremaschi, Giovannini and Schlenker (\cite[Theorem 1.4, Corollary 1.5]{CremaschiGiovanniniSchlenker}). Our results provide more examples of Schottky uniformizations with negative renormalized volume for small genus, while they give more examples for larger genus.
   
\begin{acknowledgements}
We thank the referee for their careful reading and valuable feedback, which helped clarify many parts of the paper and identify several inaccuracies. The first author was supported by the Luxembourg National research Fund AFR/Bilateral-ReSurface 22/17145118 for the latest editions of this note. The second author was partially supported by NSF grant DMS-2001997.
\end{acknowledgements}

\section{Background}

\subsection{Surfaces}
\label{subsec:surfaces}

\begin{definition}[Surface]
By a \emph{closed surface} $S$ we mean a smooth, oriented, compact 2-manifold of genus $g \geq 2$. We might also consider \emph{surfaces with boundary}, where the boundary of $S$, denoted by $\partial S$, consists of a finite union of circles. If we want to be explicit about the genus of $S$, we will write $S_g$. We will also allow non-compact surfaces with punctures in $S \backslash \partial S$.
\end{definition}

We will refer to $S$ together with a choice of metric, as $\Sigma$. If we pick a Riemannian metric and want to be explicit about the Riemannian tensor, we will write $(\Sigma,g)$.
In all the above types of topological surfaces, we will assume $S$ admits a hyperbolic metric.

\begin{definition}[Hyperbolic surface, hyperbolic subsurface]
Given a closed surface $S$, we can endow it with a choice of hyperbolic metric $\Sigma$, i.e., a complete Riemannian metric of constant curvature -1. In this case, we will call $\Sigma$ a \emph{closed hyperbolic surface}. If $S$ has non-empty boundary, we consider the boundary to be totally geodesic. If $S$ has punctures, we consider them as cusps. In this case we will say $\Sigma$ is a \emph{surface of finite type}.
Consider non-null homotopic and non-peripheral simple curves $\{ \gamma_i \}_{i=1,\cdots,n}$ on $S$ and endow $S$ with a hyperbolic metric $(\Sigma,d)$. The subspace $S' \coloneqq S - \{ \gamma_i \}_{i=1,\cdots,n} \subset S$ is a possibly disconnected open surface endowed with the restriction metric of $(\Sigma,d)$. Each connected component $S_i$ of $S'$ can be homotoped to a subsurface $S_i''$ bounded by geodesics homotopic to a subcollection of $\{ \gamma_i \}_{i=1,\cdots,n}$. The metric completion of $(S_i'',d|_{S_i''})$ is a hyperbolic metric $d'$ on $S_i$ with geodesic boundary. 
Accordingly, a \emph{hyperbolic subsurface} $\Sigma'$ of $\Sigma$ is a topological surface $S'$ so that $S' \subset S$  is endowed with complete hyperbolic metric $\Sigma'$ coming from the completion of the restriction of the metric $\Sigma$ to $\Sigma'$, as described in the previous paragraph.
\label{def:subsurface}
\end{definition}

\begin{definition}[Moduli space] Let $\mathcal{M}(S)$ denote the \emph{moduli space} of all possible hyperbolic metrics on $S$ up to oriented isometry.
\end{definition}

\subsection{Curves on surfaces}
\label{subsec:curves}

\begin{definition}[Curves]
Let $\gamma$ in a surface $S$ be a smooth connected
1-manifold $X(\gamma)$, possibly with boundary, together with a map
(also called $\gamma$) from $X(\gamma)$
into~$S$, sending $\partial X(\gamma)$ to $\partial S$. If $X(\gamma)$ is a circle, $\gamma$ will be called a \emph{curve}. If $X(\gamma)$ is a closed interval, $\gamma$ will be called an arc.
We say that $\gamma$ is \emph{non-essential} if it is contained in a disk (so it is \emph{null-homotopic}), it is \emph{boundary parallel}, i.e. homotopic to a component of $\partial S$, or if it is \emph{peripheral}, i.e., it can be homotoped into a puncture.
We say it is \emph{essential} otherwise.
Two curves $\gamma$ and~$\gamma'$ are considered equivalent if they are related by 
\emph{homotopy} within the space of all maps $X(\gamma)$ to~$S$ (not necessarily immersions), \emph{reparametrization} of the 1-manifold and
dropping trivial components.
A curve~$\gamma$ is \emph{simple} if it has a representative for which the map $\gamma$ is
injective. The space of curves on $S$ will be denoted by $\mathcal{C}(S)$.
\end{definition}

For the following lemmas, see \cite[Theorem~1.6.6]{Bus10:GeometryRiemannSurfaces} and \cite[Theorem~1.5.3]{Bus10:GeometryRiemannSurfaces}.
\begin{lemma}[Realizing curves by geodesics]
Let $\Sigma$ be a hyperbolic surface with totally geodesic boundary and punctures. If a curve $\gamma$ is essential on $\Sigma$, then $\gamma$ is represented by a unique closed geodesic $\gamma_*$, with the following properties:
\begin{enumerate}
    \item $\gamma_*$ is either contained in $\partial \Sigma$ or $\gamma_* \cap \partial \Sigma= \emptyset$
    \item If $\gamma$ is simple then $\gamma_*$ is simple
\end{enumerate}
\label{lem:georep}
\end{lemma}

\begin{lemma}\label{lem:hypintersect}
Let $\Sigma$ be a hyperbolic surface with totally geodesic boundary $\partial \Sigma$. Then any arc $\eta$ has a unique arc representative $\eta_*$ which is geodesic with respect to the hyperbolic metric and is orthogonal to $\partial \Sigma$ at the endpoints. Moreover, the geometric intersection between $\eta$ and any simple curve $\gamma$ is realized by the intersection between $\eta_*$ and the unique geodesic representative $\gamma_*$, existing by Lemma~\ref{lem:georep}.
\end{lemma}

The following Lemma is standard and follows by classification of surfaces (see~\cite[Page~27]{FB12:Primer}). 

\begin{lemma}[Existence of simple curves]
Given a hyperbolic surface of finite type which is not the pair of pants, there exists at least one essential simple closed geodesic.
\label{thm:existgeodesic}
\end{lemma}

\subsection{Hyperbolic and extremal length}
\label{subsec:ELandL}

\begin{definition}[Length of a curve]
Given a metric $\Sigma$, we define the \emph{length} of a curve $\gamma$ with respect to $\Sigma$ as
\[
\ell(\gamma;\Sigma) \coloneqq \inf_{\eta \sim \gamma} \mathrm{len}(\eta;\Sigma).
\]
where the infimum runs over curves $\eta$ equivalent to $\gamma$, and $\mathrm{len}(\gamma;\Sigma)$ denotes the length of the path $\gamma$ with respect to the metric $\Sigma$.
If $\Sigma$ is a hyperbolic metric, we call it \emph{hyperbolic length} of $\gamma$. When the underlying hyperbolic structure is understood, we will denote it just by $\ell(\gamma)$.
\end{definition}

In most of the paper we will be working with hyperbolic length $\ell(\gamma)$, but length with respect to other metrics will be relevant when referring to extremal length, which we now define.

\begin{definition}[Conformal class of metrics]
Let $\rho \colon S \to \mathbb{R}_{\geq 0}$ be a measurable function.
For a fixed $\Sigma \in \mathcal{M}(S)$, we denote by $\rho \Sigma$ the new metric on $S$ obtained by scaling the Riemannian metric tensor of $\Sigma$ by $\rho^2$.
Let $[\Sigma]$ denote the set of all such scaled metrics and refer to it as the \emph{conformal class of metrics of $\Sigma$}.
\end{definition}

We can also define a notion of length of curves which only depends on $[\Sigma]$.

\begin{definition}[Extremal length]
We define the \emph{extremal length} of a curve $\gamma$ with respect to $[\Sigma]$, for $\Sigma$ a hyperbolic metric on $S$ as

\[
\EL(\gamma;\Sigma) \coloneqq \sup_{\rho} \frac{\ell(\gamma;\rho\Sigma)^2}{\Area(\rho\Sigma)},
\]
\label{def:EL}
where the supremum is taken over all metrics $\rho\Sigma \in [\Sigma]$, so that $\ell(\gamma;\rho\Sigma)$ denotes the length with respect to the scaled metric $\rho\Sigma$ and $0<\Area(\rho\Sigma)< \infty$, where $\Area(\rho\Sigma)$ denotes the area of $S$ with respect to $\rho\Sigma$.
\end{definition}
This definition extends naturally to (weighted) multi-curves by extending the length 
$\ell(\gamma;\rho\Sigma)$ linearly and taking the same supremum as above. 
The resulting extremal length, however, is not linear as a function of the 
multi-curve.
The metric achieving the supremum in Definition~\ref{def:EL} will be called \emph{extremal metric} for $\gamma$ with respect to $\Sigma$. Although it is always achieved (see, for example, \cite[Theorem~12]{Rodin74:ExtremalLength}), the extremal metric is not well-understood for \emph{non-simple} multi-curves on a Riemann surface. However, the extremal metric has been characterized for weighted simple multi-curves.

Extremal length in general is difficult to compute, and only known explicitly in very few examples of surfaces of low complexity. Such an example is the annulus.

\begin{proposition}
If $\Sigma=A$ is a hyperbolic annulus, and $\gamma$ is in the homotopy class $C$ of the core curve of $A$, the supremum in Definition~\ref{def:EL} is realized by the Euclidean metric on the annulus. Furthermore, if the annulus is realized by identifying the horizontal edges of a Euclidean rectangle of width $a$ and height $b$ by a Euclidean translation, the extremal length of $\gamma$ in $A$ is given by

\[
  \EL(\gamma;A)=\frac{a}{b}
\]
\label{def:ELannulus}
\end{proposition}
\begin{proof}
See \cite[4.2~Examples]{Ahl10:ConformalInvariants}.
\end{proof}

By the work \cite{Jenkins57:OnExistenceOfExtremalLength}, for a weighted simple (multi-)curve (i.e., a finite union of disjoint simple curves) $\gamma$ the supremum in Definition~\ref{def:EL} is achieved by a metric which is Riemannian except at finitely many points. This metric is induced by a holomorphic (in the case of a surface with punctures, meromorphic) quadratic differential $q_C$. This object induces a decomposition of the surface into a family of Euclidean cylinders whose heights are given by the weights of the multi-curve, and its core curves are in the homotopy class of the connected components of the multi-curve (see \cite[Theorem~21.1]{Str84:QuadraticDiff} for a precise statement).
 In that case, the extremal length can can also be characterized as the $L^1$-area of the associated quadratic differential~ (see \cite[Proposition~3.12]{KPT17:ConformalEmbeddings}). In this paper we will be concerned with the extremal length of connected curves.

\subsection{Previous comparison results}
\label{subsec:previouscomparison}

For a fixed point $\Sigma \in \mathcal{M}(S)$, and an arbitrary curve $\gamma \in \mathcal{C}(S)$, the relation between $\ell(\gamma;\Sigma)$ and $\EL(\gamma;\Sigma)$ is complicated.
Comparison results between extremal length and hyperbolic length have been obtained before. Maskit \cite[Propositions~1~and~2]{Maskit85:ComparisonLengths} obtains lower and upper bounds for arbitrary conformal classes and curves.

\begin{theorem}[Maskit comparison result]\label{thm:maskit}
For any $\gamma \in \mathcal{C}(S)$, and any $\Sigma \in \mathcal{M}(S)$,
\[
\frac{\ell(\gamma;\Sigma)}{\pi}  \leq \EL(\gamma;\Sigma) \leq \frac{1}{2}\ell(\gamma;\Sigma)e^{\ell(\gamma;\Sigma)/2}.
\]
\label{thm:MaskitComparison}
\end{theorem} 
Thus, for a fixed curve $\gamma$, if $\ell(\gamma;\Sigma) \to 0$ as a function of moduli space, it follows that $\EL(\gamma;\Sigma) \to 0$, too, and in fact $\lim_{\ell \to 0} \frac{\EL}{\ell} = \frac1\pi$.

Maskit's upper estimate is close to being sharp: in \cite{Maskit85:ComparisonLengths}, Maskit produces a family of examples of uniformized punctured tori for which the extremal length of the free homotopy class of one of the standard generators of the fundamental group is lower bounded by an exponential function of the hyperbolic length. This upper bound follows by Maskit's estimate, contained in the proof of Theorem~\ref{thm:upperEL}.

On the other hand, Minsky \cite{Min96:ExtremalLength} gives estimates of the extremal length of a curve on a Riemann surface in terms of the contributions of the hyperbolic lengths of the restriction of the curve to  subsurfaces. The result shows that the contributions are of two types:

\begin{enumerate}
\item Of order the square of the hyperbolic length in the regions with large injectivity radius (these regions are pairs of pants).
\item Depending on the modulus and amount of twisting of the curve around the annular regions with small injectivity radius. In the case where the curve is just homotopic to the core curve of the annuli, the contribution is just proportional to the inverse of the modulus which is also proportional to the hyperbolic length.
Note how this is compatible with Maskit's result discussed above.
\end{enumerate}

\subsection{Hyperbolic and extremal length systoles}
\label{subsec:systoles}
In this section we motivate Theorem~\ref{thm:mainthm} by focusing on the first curve in a short partial pants decomposition, the hyperbolic systole. We show that for that curve one can obtain a universal lower bound on the width of an annular collar directly applying the work of Buser and Sarnak.  We also introduce the concept of extremal length systole and make some initial observations to compare it to its hyperbolic analogue.

As we mentioned in the introduction, and explain precisely in Definition~\ref{def:shortdecomposition}, the first curve of a short partial pants decomposition is one of the (possibly multiple) shortest curves with respect to the hyperbolic metric $\Sigma$, whose length we call the \emph{hyperbolic length systole}. Precisely, we define
\[
\sys_{\ell}(\Sigma) \coloneqq \inf_{\gamma \in \mathcal{S}(S)} \ell(\gamma;\Sigma),
\]
where the infimum is taken over $\mathcal{S}(S)$, all essential simple curves in $S$.
We can define the analogous quantity for extremal length, this is the \emph{extremal length systole}
\[
\sys_{\EL}(\Sigma) \coloneqq \inf_{\gamma \in \mathcal{S}(S)} \EL(\gamma;\Sigma).
\]

Although these two quantities are different, they share some parallelisms. 

For example, we note that for any metric in $[\Sigma]$, as long as $\Sigma$ is not the thrice punctured sphere or a pair of pants, resolving an essential crossing on a geodesic curve $\gamma$ can never increase the length of the geodesic representative of the resulting curve. Therefore, taking the above infima over all curves or just over all simple curves yields the same extremal and hyperbolic length systoles (see \cite[Lemma~4.16]{MGT20:FromCurvesToCurrents}).

The main goal of this paper is to obtain a uniform linear bound for the extremal length in terms of the hyperbolic length of curves of a short pair of pants decomposition. 
We observe that the comparison results outlined in Subsection~\ref{subsec:previouscomparison} are not enough to obtain the desired uniform linear bound, but at the same time they do not use the fact that the curves are short.
We will say that two functions of the genus, $\phi,\psi : \mathbb{N} \to \mathbb{R}$ are \emph{asymptotic}, and write $\psi \asymp \phi$ if there exist positive constants $C_2>C_1$ independent of the genus such that, for all large enough $g$,
\[
      C_1 \psi(g) \leq \phi(g) \leq C_2 \psi(g).
\]

By results of Gromov (for upper bound) \cite[Theorem~5.5~C']{Gro83:FillingManifolds} and its proof, and Buser-Sarnak (for lower bound) \cite[Section~4]{BS94:PeriodMatrix}, the following asymptotics for the hyperbolic length systole are known.
\begin{theorem}[Gromov, Buser-Sarnak asymptotics for hyperbolic length systole]
\[
\max_{\Sigma \in \mathcal{M}(S)} \sys_{\ell}(\Sigma) \asymp \log(g).
\]
\label{thm:sysllog}
\end{theorem}
If $\gamma$ is a simple closed geodesic in $\Sigma$, for a given $w>0$, we consider the set
\[
\widehat{A_w(\gamma)} \coloneqq \{ x \in \Sigma : d(x,\gamma) \leq w \}.
\]
For $w$ small enough, $\widehat{A_w(\gamma)}$ is an embedded annulus. Choosing $w$ to be smaller or equal than the largest number with this property, we will call $\widehat{A_w(\gamma)}$ a \emph{annular collar} (or just \emph{collar}).
Note that $\gamma$ divides $A_{w}(\gamma)$ into two annuli that we will refer to as \emph{half collars}, and denote by $A_w(\gamma)$.

The next result follows essentially by work of Buser and Sarnak in \cite[(3.12)]{BS94:PeriodMatrix}, which we will present for completeness in the beginning of proof of Lemma~\ref{lem:ortho}.

\begin{theorem}[Buser-Sarnak, 94]
Let $\Sigma$ be a closed hyperbolic surface.
Let $\gamma$ be the hyperbolic systole of $\Sigma$. There exists a universal constant $w$ independent of the genus $g$, so that there is an annular collar $\widehat{A_w(\gamma)}$.
\label{thm:busersarnak}
\end{theorem}

As an immediate corollary of this result and Theorem~\ref{thm:upperEL}, we observe that there is a linear upper bound of extremal length of the hyperbolic systole, and thus of the extremal length systole, in terms of its hyperbolic length.
From this and the lower bound of Theorem~\ref{thm:maskit}, it follows that extremal length systole and hyperbolic length systole have the same asymptotic behaviour for large genus.

\begin{corollary}
Let $\Sigma$ be a closed hyperbolic surface. Then
\[
\max_{\Sigma \in \mathcal{M}(S)} \sys_{\EL}(\Sigma)  \asymp \log(g).
\]
\label{cor:sysellog}
\end{corollary}

\subsection{Renormalized volume}
\label{subsec:renormalizedvolume}
One of the main applications of our comparison result is Corollary~\ref{cor:vrsystole} which gives upper bounds of the renormalized volume of certain 3-manifolds in terms of the hyperbolic length of a family of compressible curves.
Renormalized volume for hyperbolic $3$-manifolds (as described in \cite{KS08}) is motivated by the computation of the gravity action $S_{gr}[g]$ in the context of the Anti-de-Sitter/Conformal Field Theory (AdS/CFT) correspondence \cite{Witten98}. Let us start by defining the $W$-$volume$ of a compact, convex $C^{1,1}$-submanifold $N\subset M$, of a convex-cocompact hyperbolic 3-manifold $M$, by
\begin{equation}
    W(M,N)= {\rm vol}(N) - \frac12 \int_{\partial N}Hda,
\end{equation}
where ${\rm vol}(N)$ is the hyperbolic volume of $N$ with respect to the metric in $M$, $H$ is the mean curvature of $\partial N$, and $da$ is the area form of the induced metric on $\partial N$.
The \emph{boundary at infinity of $M=\mathbb{H}^3/\Gamma$}, with $\Gamma = \pi_1(M)$, is defined as the surface
$\partial_\infty M \coloneqq \Omega(\Gamma)/\Gamma$, where $\Omega(\Gamma)=\partial \mathbb{H}^3
 \backslash\Lambda(\Gamma)$ is the domain of discontinuity of $\Gamma$,
and $\Lambda(\Gamma)$ the limit set of the action. Since $\Omega(\Gamma)$ is an open subset of $\partial \mathbb{H}^3$ and $\PSL(2,\mathbb{C})$ acts by biholomorphisms of $\partial \mathbb{H}^3$, it follows that $\partial_\infty M$ is 
equipped with a natural complex projective structure, and hence
a complex structure or---equivalently by the Uniformization Theorem---a conformal class of metrics on $\partial_\infty M$, which by abuse of notation we also call $\partial_\infty M$.
Given a metric $h$ in the conformal class $\partial_\infty M$, Epstein \cite{Epstein} constructs a family of homotopically equivalent
convex submanifolds $N_r$ with equidistant boundary by taking envelopes of horospheres. Such family depends on the projective structure of $\partial_\infty M$ and exhausts the ends of $M$. Because the boundaries are equidistant, the $W$-volumes have the property \cite[Lemma 3.6]{Schlenker13}
\begin{equation}
    W(M,N_r) = W(M,N_s) - \pi(r-s)\chi(\partial M),
\end{equation}
which leads to the definition (independent of $r$):
\begin{equation}
    W(M,h) := W(M,N_r) + \pi r\chi(\partial M).
\end{equation}
Taking $h_{\rm hyp}$ the metric of constant curvature in the given conformal class at infinity, we define \emph{Renormalized Volume} $\VR$ as
\begin{equation}
    \VR(M) \coloneqq W(M,h_{\rm hyp}).
\end{equation}
As an example, one can take $M=M_0\cup E_1\cup\ldots\cup E_k$, where $M_0$ is a compact hyperbolic $3$-manifold with totally geodesic boundary $\partial M_0 = (\Sigma_1,g_1)\cup\ldots\cup(\Sigma_k,g_k)$ and $\lbrace E_i = \Sigma_i\times\mathbb{R}^{\geq0}\rbrace_{1\leq i \leq k}$ are the ends of $M$ attached to $M_0$ through its boundary and with geometry $(\Sigma_i\times\mathbb{R}^{\geq0}, \cosh^2(t)g_i + dt^2 )$. In this particular case the product structure in each end $E_i=\Sigma_i\times\mathbb{R}^{\geq 0}$ is the equidistant foliation associated to $h_{\rm hyp}$, so that $N_0$ corresponds to $M_0$. Therefore, we calculate $\VR(M)$ as
\begin{equation}
V_R(M) = W(M,N_0) = {\rm vol}(M_0) - \frac12 \int_{\partial M_0}Hda = {\rm vol}(M_0)
\end{equation}
When $M$ is convex co-compact, its compactification is given by $M \cup \partial_\infty M$, where $\partial_\infty M$ is homeomorphic to the closed surface $S=\partial \widehat{M}$. Hence the conformal boundary determines a point $[\partial_\infty M] \in \Teich(S)$.

Naturally, $V_R$ is a function on $\operatorname{CC}(M)$, the space of convex co-compact hyperbolic structures on $M$ up to homotopy equivalence. By~\cite[Theorem~5.1.3]{Marden2016:HyperbolicManifolds} and~\cite[Theorem~5.27]{MatsuzakiTaniguchi1998:HyperbolicManifolds} it follows that $\operatorname{CC}(M)$ is parameterized by $\Teich(\partial \overline{M})/\sim$, where $\Teich(\partial \overline{M})$ is the product of the Teichm\"uller spaces of the components of $\partial \overline{M}$, and the equivalence relation $\sim$ is generated by Dehn twists along compressible curves in $\partial \overline{M}$.
$V_R$ is thus naturally a function on Teichm\"uller space of the boundary, and, in fact, can be used to obtain a geometric K\"ahler potential for the Weil-Petersson metric in Teichm\"uller space, as well as an entropy for hyperbolic 3-manifold with a given conformal boundary (see for instance \cite{KS08}). For this reason, estimates for renormalized volume are of interest and can be used to test intuition coming from theoretical physics.
In this paper, we will use our comparison result to address the following question, due to Maldacena.

\begin{question}[Maldacena]\label{que:maldacena}
For which Riemann surfaces $(\Sigma,g)$ can we find a handlebody $M^3$ so that $\Sigma$ is conformal to the conformal class at infinity of $M$ and
\[2\VR(M) \leq \VR(\Sigma\times \mathbb{R},\cosh^2(t)g + dt^2).
\]
\end{question}
The motivation for this question comes from the intuition that the pair of handlebodies should be more `likely' than a Fuchsian filling of two copies of $\Sigma$ by the product, where in more precise terms `likehood' refers to the term in the partition function with the largest contribution, which is given by the largest value for $e^{-\VR}$~\cite{Krasnov2000,SchlenkerWitten2022}. Note that from our example, we can find that $\VR(\Sigma\times \mathbb{R},\cosh^2(t)g + dt^2) = 0$. In Corollary \ref{cor:vrsystole} we will see that the answer is yes as long as $\Sigma$ has a suitable collection of curves (see paragraph above Theorem~\ref{mainthm:renormalized}), i.e., sufficiently many short curves decomposing the surface into subsurfaces of either genus 1, or genus 2 with a boundary component.

For a more detailed description of renormalized volume for hyperbolic $3$-manifolds, we refer the reader to \cite{KS08}. For other bound results on $\VR$ one has the following references \cite{Schlenker13}, \cite{BC15}, \cite{BBB19}, \cite{VP17}, \cite{VP19:UpperBounds},\cite{CremaschiGiovanniniSchlenker}.

\subsection{Main results}

The main results of this paper are the following.
Recall that an $\epsilon$-short partial pants decomposition is a collection of $n < 3g-3$
disjoint simple closed geodesics with $\ell(\gamma_i)<\epsilon$.
The map $f$ is defined in the proof of Lemma~\ref{lem:noshortorthogeodesics}.

\begin{theorem}\label{thm:mainthm}
There exist universal constants
\[
w_0\simeq 0.40352, \qquad
c = \tfrac{2}{\pi-2\theta} \simeq \tfrac{2}{\pi-2\cdot 1.1778} \simeq \tfrac{2}{\pi-2.3556}, \qquad
\epsilon=f^{-1}(w_0)\approx 3.2282
\]
with the following property.  
Let $\Sigma$ be a closed hyperbolic surface, and let
$\mathcal{C}=\{\gamma_i\}_{1\leq i \leq n}$, for $n<3g-3$, be an $\epsilon$-short partial pants decomposition. Let $\gamma$ be
the shortest essential curve in a component of
$\Sigma \setminus \mathcal{C}$ that contains at least one simple essential
curve. Then $\gamma$ admits an embedded annular half collar
$A_{w_0}(\gamma)$ of width $w_0$. In particular, its extremal length satisfies
\[
  \EL(\gamma;\Sigma) \leq c\,\ell(\gamma;\Sigma).
\]
\end{theorem}

We also obtain upper bounds for the extremal length of hyperbolic systole of a surface with boundary.

\begin{proposition}
Let $\Sigma$ be a hyperbolic surface with totally geodesic boundary.
Assume that each component of $\partial\Sigma$ has a collar of width $w$, and all such collar are disjoint. Let $\gamma$ be one of the shortest essential hyperbolic geodesics in $\Sigma$. Then 
\begin{equation}
\EL(\gamma;\Sigma)  \leq \frac{4\ell(\gamma)}{\pi-2\theta}\
\end{equation}
where $\theta\in[0,\pi/2)$ satisfies $\sin(\theta)\cdot\cosh(\min\lbrace w_0,w\rbrace)=1$ with $w_0\simeq 0.40352$ 
\label{prop:boundary_systole}
\end{proposition}

Finally, we comment that the question of what the global maximum of the hyperbolic systole is for a specific genus $g$ has also been considered in the literature. For example, for genus $g=2$, $\max \sys_{\ell}$ is realized by the Bolza surface (see \cite[Theorem~5.2]{Sch93:ShortestGeodesic}, for example).
The Bolza surface also appears as a maximum for systoles of other notions of lengths (see, for example \cite[Appendix~A.1]{BS94:PeriodMatrix}).

We anticipate the following additional parallelism between the extremal and hyperbolic length systoles.

\begin{conjecture}
For genus 2, $\max \sys_{\EL}$ is attained by the conformal class of the Bolza surface.
\end{conjecture}

In work with Maxime Fortier Bourque \cite{FBMGVP:Bolza} we showed the following.

\begin{theorem}
For genus 2, $\sys_{\EL}$ attains a local maximum at the conformal class of the Bolza surface.
\end{theorem}

\section{Proof of Main theorem}
\label{subsec:results}

In this section we extend Buser-Sarnak techniques in \cite[Section~4]{BS94:PeriodMatrix} to give a linear upper bound of the extremal length of curves in certain short partial pair of pants decomposition in terms of their hyperbolic length.

\begin{definition}[Short partial pair of pants decomposition]
Let $\Sigma$ be a closed hyperbolic surface.
A \emph{short partial pair of pants} $\Gamma \coloneqq \lbrace \gamma_i\rbrace_{1\leq i\leq n}$, where $n\leq 3g-3$ is a set of disjoint simple curves in a Riemann surface $\Sigma$ obtained inductively as follows. Given $\lbrace \gamma_i\rbrace_{1\leq i\leq k-1}$ we define $\gamma_k$ as (one of possibly many) shortest closed geodesic in one component of the metric completion (as in Definition~\ref{def:subsurface}) of $\Sigma\setminus \lbrace \gamma_i\rbrace_{1\leq i\leq k-1}$ which is essential.
\label{def:shortdecomposition}
\end{definition}

\begin{lemma}\label{lem:noshortorthogeodesics}
Given two disjoint curves $\delta,\gamma$ in a surface $\Sigma$ so that $\gamma$ is simple, $\epsilon=\ell(\delta)$, there exists a function $f(\epsilon)$, so that any orthogeodesic $\eta$ between $\delta$ and $\gamma$ has length greater than $f(\epsilon)$. In particular, if $\gamma$ is at distance less than $f(\epsilon)$ from $\delta$, then $\gamma$ and $\delta$ must intersect.

\end{lemma}

\begin{proof}
Consider the upper half plane model of $\mathbb{H}^2$. Lift to $\mathbb{H}^2$ so that $\eta$ corresponds to the unit circle arc of angle $\theta$ between $(0,1)$ and $(\sin\theta,\cos\theta)$ for $0<\theta<\pi/2$, the vertical $y$-axis quotients to $\delta$ by the map $z\mapsto e^\epsilon z$, and the upper half-circle $\tilde{\gamma}$ orthogonal to the unit circle at $(\sin\theta,\cos\theta)$ is an axis for $\gamma$. From simple calculations it follows that the length of $\eta$ is expressed as $\ln(\sec\theta+\tan\theta)$ and the first coordinates of the  axis $\tilde{\gamma}$ are $\csc\theta-\cot\theta$, $\csc\theta+\cot\theta$. Since $\gamma$ is simple, we know that $\tilde{\gamma}\cap e^\epsilon\tilde{\gamma}=\emptyset$, from where it follows then that $\frac{\csc\theta+\cot\theta}{\csc\theta-\cot\theta} = \left( \csc\theta + \cot\theta\right)^2<e^\epsilon$. As for $0<\theta<\pi/2$ we have that $g_1(\theta)=\csc\theta+\cot\theta$ is strictly decreasing and $g_2(\theta):=\sec\theta+\tan\theta$ is strictly increasing, then the length of $\eta$ is bounded by
\[\ell(\eta) > \ln \left(g_2\left(g^{-1}_1(e^{\epsilon/2})\right)\right) =:f(\epsilon).
\]
\end{proof}

We proceed now with the proof of Theorem~\ref{thm:mainthm}.

Recall that a simple curve $\gamma$ is \emph{separating} if $S\backslash \gamma$ is a disconnected surface. If a curve is separating then, for any curve $\eta$ that intersects $\gamma$, the intersection number is at least $2$.

\begin{lemma}\label{lem:ortho}
Let $\Sigma$ be a compact hyperbolic surface with totally geodesic boundary so that $\Sigma \neq \Sigma_{0,3}$. Let $\gamma$ be one of the simple closed geodesics of least length among the set of essential ones. Let $\eta$ be an orthogeodesic with endpoints in $\gamma$. 

\begin{enumerate}
\item $\Sigma \neq S_{0,4}$.
\begin{itemize}
    \item If $\eta$ arrives at opposite sides of $\gamma$ (and in particular $\gamma$ is non-separating), then $\ell(\eta)\geq w_1=2\arctanh(\frac23)\approx 1.6094$
    \item If $\eta$ arrives at the same side of a non-separating curve $\gamma$, then $\ell(\eta)\geq w_2 \coloneqq 2\arcsinh(\sqrt{3}) \approx 2.634$
    \item If $\eta$ arrives at the same side and $\gamma$ is separating, then there is at least one component $\Sigma^-$ of $\Sigma\setminus\gamma$ which is necessarily not a pair of pants and $\ell(\eta) \geq 2\arcsinh(\sqrt{3}) \approx 2.634$ for $\eta$  in $\Sigma^-$. 
\end{itemize}
\item $\Sigma=S_{0,4}$. Here $\eta_{1,2}$ are the shortest orthogeodesics on each side, and $\gamma$ is separating.
\begin{itemize}
\item If the endpoints $\eta_{1,2}$ alternate:
$\max \{\ell(\eta_1),\ell(\eta_2)\} \geq w_3 \simeq 0.80507$
\item If the endpoints $\eta_{1,2}$ do not alternate:
$\max \{\ell(\eta_1),\ell(\eta_2)\} \geq w_4 \simeq 1.2188$.
\end{itemize}
\end{enumerate}
\end{lemma}

\begin{proof}
Since $\Sigma \neq \Sigma_{0,3}$, the set of essential simple closed geodesics is non-empty (by Lemma~\ref{thm:existgeodesic}). Let $\gamma$ be the shortest one. As in \cite{BS94:PeriodMatrix}, we can use a short orthogeodesic $\eta$ to produce a curve $\tilde{\gamma}$ in $\Sigma$ which is strictly shorter than $\gamma$. The challenge is to prove, once we have defined $\tilde{\gamma}\subset\Sigma$, that it is essential. We analyze this depending on if $\eta$ arrives at opposite/same sides of $\gamma$, the former being possible only when $\gamma$ is non-separating. We will deal with the case $\Sigma = S_{0,4}$ independently. We observe that the methods from \cite{BS94:PeriodMatrix} extend to this case, but we include their arguments for sake of completeness.

\begin{itemize}
    \item $\eta$ arrives at opposite sides (and thus $\gamma$ is non-separating). Take two lifts of $\gamma$ in $\mathbb{H}^2$ joined by the geodesic $\eta$ orthogonal to both of them. Let $\delta$ be the geodesic joining points in each lift of $\gamma$ that are equidistant to $\eta$ and correspond to the same point in the quotient (see Figure \ref{fig:lifts}). The distance of any of these lifts of points to $\eta$---which could be unique--- is at most $\frac14\ell(\gamma)$. Applying the cosine formula (see Theorem~\ref{eq:coslaw}) we have
    
    \begin{equation}
        \cosh\frac14\ell(\gamma)\cdot\cosh\frac12 \ell(\eta) \geq \cosh\frac12\ell(\delta).
    \end{equation}
    Since $\delta$ is a geodesic that meets the lifts of $\gamma$ at the same angles, its quotient is a closed geodesic in $\Sigma$ with intersection number $1$ with $\gamma$. Hence $\delta$ is essential, and the inequality $\ell(\delta)\geq\ell(\gamma)$ holds.  In particular we get
    
    \begin{equation}
        \cosh\frac14\ell(\gamma)\cdot\cosh \frac12\ell(\eta) \geq \cosh\frac{\ell(\gamma)}{2}.
    \end{equation}
    
    On the other hand, we also have, by the Collar's Lemma~\cite[4.1.1.~Theorem]{Bus10:GeometryRiemannSurfaces} that $\sinh\frac12 \ell(\gamma)\cdot\sinh \frac12\ell(\eta) > 1$. Both inequalities then imply
    
    \begin{equation}
        \frac12 \ell(\eta) \geq \max \Bigg\lbrace \arccosh\left( \frac{\cosh\frac12\ell(\gamma)}{\cosh\frac14\ell(\gamma)}\right),  \arcsinh\left( \frac{1}{\sinh\frac12\ell(\gamma)}\right) \Bigg\rbrace.
    \end{equation}

    Since the first function in the
maximum is increasing in $\ell$ while the second is decreasing, the minimum occurs at
the point where the two functions assume the same value, which is $\arctanh(\frac23)$. Hence we have that $\frac12 \ell(\eta)$ is bounded below by $\arctanh(\frac23)$, or equivalently, $\ell(\eta)\geq w_1 := 2\arctanh(\frac23)$.

    \begin{figure}[H]
\centering{
\resizebox{80mm}{!}{\Huge{
\begingroup%
  \makeatletter%
  \providecommand\color[2][]{%
    \errmessage{(Inkscape) Color is used for the text in Inkscape, but the package 'color.sty' is not loaded}%
    \renewcommand\color[2][]{}%
  }%
  \providecommand\transparent[1]{%
    \errmessage{(Inkscape) Transparency is used (non-zero) for the text in Inkscape, but the package 'transparent.sty' is not loaded}%
    \renewcommand\transparent[1]{}%
  }%
  \providecommand\rotatebox[2]{#2}%
  \newcommand*\fsize{\dimexpr\f@size pt\relax}%
  \newcommand*\lineheight[1]{\fontsize{\fsize}{#1\fsize}\selectfont}%
  \ifx\svgwidth\undefined%
    \setlength{\unitlength}{339.26319753bp}%
    \ifx\svgscale\undefined%
      \relax%
    \else%
      \setlength{\unitlength}{\unitlength * \real{\svgscale}}%
    \fi%
  \else%
    \setlength{\unitlength}{\svgwidth}%
  \fi%
  \global\let\svgwidth\undefined%
  \global\let\svgscale\undefined%
  \makeatother%
  \begin{picture}(1,1)%
    \lineheight{1}%
    \setlength\tabcolsep{0pt}%
    \put(0,0){\includegraphics[width=\unitlength,page=1]{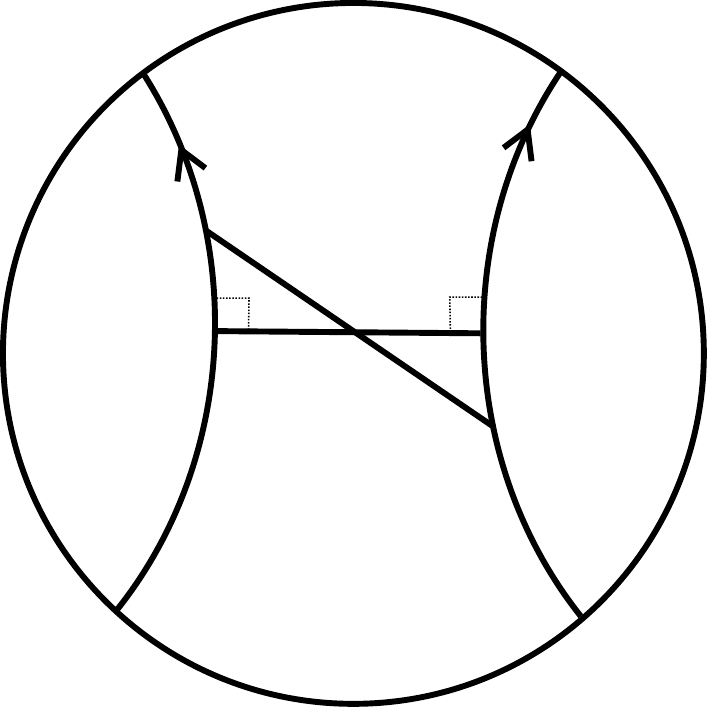}}%
    \put(0.19178202,0.70300241){\color[rgb]{0,0,0}\makebox(0,0)[lt]{\lineheight{1.25}\smash{\begin{tabular}[t]{l}$\gamma$\end{tabular}}}}%
    \put(0.72452525,0.71334238){\color[rgb]{0,0,0}\makebox(0,0)[lt]{\lineheight{1.25}\smash{\begin{tabular}[t]{l}$\gamma$\end{tabular}}}}%
    \put(0.55369172,0.42192051){\color[rgb]{0,0,0}\makebox(0,0)[lt]{\lineheight{1.25}\smash{\begin{tabular}[t]{l}$\delta$\end{tabular}}}}%
    \put(0.51263823,0.54477915){\color[rgb]{0,0,0}\makebox(0,0)[lt]{\lineheight{1.25}\smash{\begin{tabular}[t]{l}$\eta$\end{tabular}}}}%
  \end{picture}%
\endgroup%
}}
\caption{Lifts $\gamma$,$\delta$ and the orthogeodesic $\eta$ to $\mathbb{H}^2$. We abuse notation and use the same names as the curves the lifts are covering.}
\label{fig:lifts}
}
\end{figure}
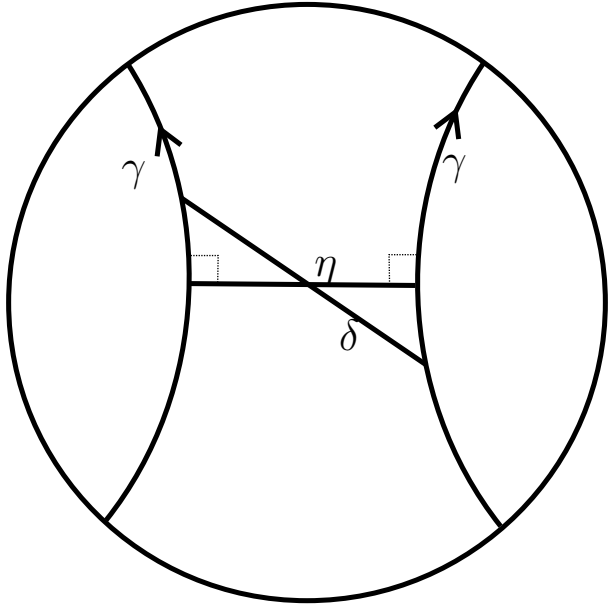

    \begin{figure}[H]
\centering{
\resizebox{80mm}{!}{\Huge{
\begingroup%
  \makeatletter%
  \providecommand\color[2][]{%
    \errmessage{(Inkscape) Color is used for the text in Inkscape, but the package 'color.sty' is not loaded}%
    \renewcommand\color[2][]{}%
  }%
  \providecommand\transparent[1]{%
    \errmessage{(Inkscape) Transparency is used (non-zero) for the text in Inkscape, but the package 'transparent.sty' is not loaded}%
    \renewcommand\transparent[1]{}%
  }%
  \providecommand\rotatebox[2]{#2}%
  \newcommand*\fsize{\dimexpr\f@size pt\relax}%
  \newcommand*\lineheight[1]{\fontsize{\fsize}{#1\fsize}\selectfont}%
  \ifx\svgwidth\undefined%
    \setlength{\unitlength}{433.31406bp}%
    \ifx\svgscale\undefined%
      \relax%
    \else%
      \setlength{\unitlength}{\unitlength * \real{\svgscale}}%
    \fi%
  \else%
    \setlength{\unitlength}{\svgwidth}%
  \fi%
  \global\let\svgwidth\undefined%
  \global\let\svgscale\undefined%
  \makeatother%
  \begin{picture}(1,0.30462419)%
    \lineheight{1}%
    \setlength\tabcolsep{0pt}%
    \put(0,0){\includegraphics[width=\unitlength,page=1]{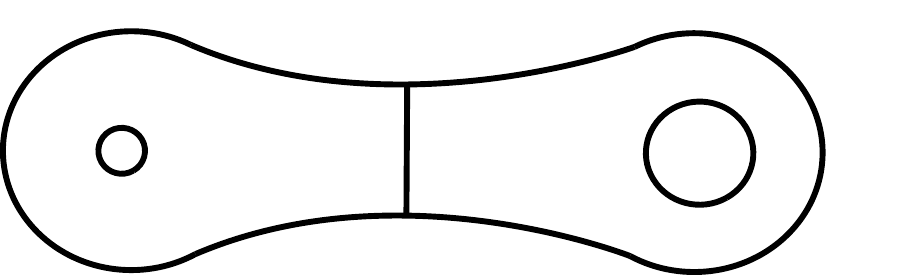}}%
    \put(0.72559126,0.28358359){\color[rgb]{0,0,0}\makebox(0,0)[lt]{\lineheight{1.25}\smash{\begin{tabular}[t]{l}$\gamma$\end{tabular}}}}%
    \put(0.45697472,0.15685225){\color[rgb]{0,0,0}\makebox(0,0)[lt]{\lineheight{1.25}\smash{\begin{tabular}[t]{l}$\eta$\end{tabular}}}}%
    \put(0,0){\includegraphics[width=\unitlength,page=2]{pantsverbose.pdf}}%
    \put(0.70823072,0.21640097){\color[rgb]{0,0,0}\makebox(0,0)[lt]{\lineheight{1.25}\smash{\begin{tabular}[t]{l}$\tilde{\gamma}^a$\end{tabular}}}}%
    \put(0.09120967,0.2032902){\color[rgb]{0,0,0}\makebox(0,0)[lt]{\lineheight{1.25}\smash{\begin{tabular}[t]{l}$\tilde{\gamma}^b$\end{tabular}}}}%
  \end{picture}%
\endgroup%
}}
\caption{Closed geodesics $\tilde{\gamma}^{a,b}$ obtained by concatenating $\gamma^a\eta, \gamma^b\eta^{-1}$. Together with $\gamma$, they are the boundary of a pair of pants in $\Sigma$ }
\label{fig:pantsverbose}
}
\end{figure}
    \item $\eta \subset \Sigma\setminus\gamma$ arrives at the same side: Let $w \coloneqq \frac12\ell(\eta)$ and denote by $\gamma^{a,b}$ the arcs in $\gamma$ between the endpoints of $\eta$. Then the concatenations $\gamma^a\eta, \gamma^b\eta^{-1}$ correspond to close geodesics $\tilde{\gamma}^{a,b}$ in $\Sigma$ (which are different in $\Sigma\setminus\gamma$ but not necessarily different in $\Sigma$) that together with $\gamma$ are the boundary of a pair of pants in $\Sigma$ (see Figure~\ref{fig:pantsverbose}).
    Assume first at least one of $\tilde{\gamma}^{a,b}$ is essential. 
    
    Without loss of generality, suppose that $\ell(\tilde{\gamma}^a)\geq \ell(\tilde{\gamma}^b)$. Since at least one of $\tilde{\gamma}^{a,b}$ is essential it follows that $\ell(\tilde{\gamma}^a)\geq \ell(\gamma)$.
    
    Denote by $t=\frac12\ell(\gamma^b), r=\frac12\ell(\gamma), r'=\frac12\ell(\tilde{\gamma}^b)$. Then by the right-angled pentagon formula~\ref{eq:rightanglepent}, we have
    
    \begin{align}
        \label{pent1}\sinh w \cdot \sinh t = \cosh r'\\
        \sinh w \cdot \sinh (r-t) = \cosh\frac12 \ell (\tilde{\gamma}^a) \geq \cosh r
    \end{align}
    which combined with $\sinh(r-t) = \sinh r\cosh t - \cosh r\sinh t$ leads to
    
    \begin{equation}
        \sinh w \cdot \cosh t \geq \coth r (1+\cosh r') \geq 1+\cosh r'
    \end{equation}
    
    Squaring this last inequality and replacing $\cosh^2 t = \sinh^2 t + 1$ we have that
    
    \begin{equation}
        \sinh^2w(\sinh^2 t + 1) \geq (1+\cosh r')^2 = 1+2\cosh r' + \cosh^2r'
    \end{equation}
    where after replacing (\ref{pent1}) we have
    
    \begin{equation}
        \sinh^2 w \geq 1+2\cosh r'\geq 3
    \end{equation}
    so then $w\geq \arcsinh\sqrt{3}$.
    
    If both $\tilde{\gamma}^{a,b}$ are non-essential then the component of $\Sigma\setminus\gamma$ that contains $\eta$ is a pair of pants and hence $\gamma$ is separating. This finishes the argument for $\gamma$ non-separating, as similarly to the previous case we have $\ell(\eta)\geq w_2:=2\arcsinh(\sqrt{3})$. For $\gamma$ separating it says that $\eta$ can only be in the unique component of $\Sigma\setminus\gamma$ that is a pair of pants, if any. So if $\Sigma\neq S_{0,4}$ the inequality $\ell(\eta)\geq w_2$ follows for at least one component of $\Sigma\setminus\gamma$.
    
    \begin{figure}[H]
\centering{
\resizebox{120mm}{!}{\Huge{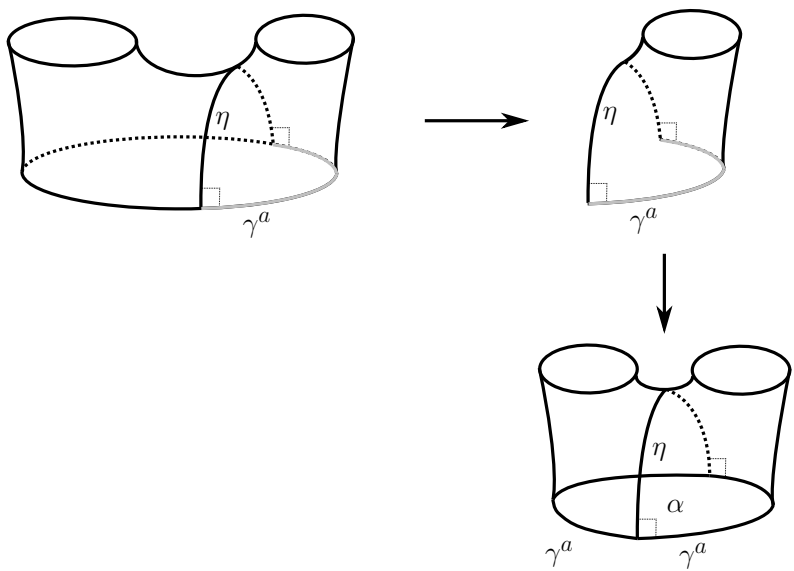}}
\caption{Doubling along $\eta_{1,2}$ to obtain a pair of pants $P$ with a cuff $\alpha$ of length $2\ell(\gamma_1)$ and orthogeodesics $\eta_{1,2}$.}
\label{fig:doubling}
}
\end{figure}
\end{itemize}

Finally, let us address the case when $\Sigma=S_{0,4}$. Notice that by the topology of $\Sigma$ there are only two simple orthogeodesics $\eta_{1,2}$ to $\gamma$, one at each side of $\Sigma\setminus\gamma$. We consider two subcases: 

\begin{itemize}
    \item The endpoints of $\eta_{1,2}$ alternate (see Figure~\ref{fig:fourholedlinked}): In this case, from the four arcs between these endpoints (and allowing for $0$ length arcs) there are two opposite arcs with lengths adding up to at most $\frac12 \ell(\gamma)$. We can concatenate them with $\eta_{1,2}$ so as to obtain a simple curve $\tilde{\gamma}$ with $\ell(\tilde{\gamma}) \leq \ell(\eta_1) + \frac12 \ell(\gamma) + \ell(\eta_2)$. Since the geometric intersection between $\gamma$ and $\tilde{\gamma}$ is $2$, then $\tilde{\gamma}$ is essential, from which we know $\ell(\gamma) \leq \ell(\tilde{\gamma})$ and subsequently $\ell(\eta_1) + \ell(\eta_2) \geq \frac12 \ell(\gamma)$. On the other hand, we can apply the Collar Lemma~\cite[4.1.1.~Theorem]{Bus10:GeometryRiemannSurfaces} for the double of $\Sigma$ to obtain half collars at each side of $\gamma$ of width $w$ satisfying $\sinh(w)\sinh(\frac12 \ell(\gamma))=1$. As orthogeodesics are obstructions for half collars, we have that $\ell(\eta_1), \ell(\eta_2) \geq 2w = 2\arcsinh(1/\sinh(\frac12\ell(\gamma)))$. Hence $\max \lbrace \ell(\eta_1), \ell(\eta_2)\rbrace \geq \max \lbrace \frac14\ell(\gamma), 2\arcsinh(1/\sinh(\frac12\ell(\gamma)))\rbrace = w_3\simeq 0.80705$, where $w_3$ is the positive constant satisfying $\sinh(2w_3)\sinh(w_3/2)=1$.
        \begin{figure}[H]
\centering{
\resizebox{80mm}{!}{\Huge{
\begingroup%
  \makeatletter%
  \providecommand\color[2][]{%
    \errmessage{(Inkscape) Color is used for the text in Inkscape, but the package 'color.sty' is not loaded}%
    \renewcommand\color[2][]{}%
  }%
  \providecommand\transparent[1]{%
    \errmessage{(Inkscape) Transparency is used (non-zero) for the text in Inkscape, but the package 'transparent.sty' is not loaded}%
    \renewcommand\transparent[1]{}%
  }%
  \providecommand\rotatebox[2]{#2}%
  \newcommand*\fsize{\dimexpr\f@size pt\relax}%
  \newcommand*\lineheight[1]{\fontsize{\fsize}{#1\fsize}\selectfont}%
  \ifx\svgwidth\undefined%
    \setlength{\unitlength}{394.75866603bp}%
    \ifx\svgscale\undefined%
      \relax%
    \else%
      \setlength{\unitlength}{\unitlength * \real{\svgscale}}%
    \fi%
  \else%
    \setlength{\unitlength}{\svgwidth}%
  \fi%
  \global\let\svgwidth\undefined%
  \global\let\svgscale\undefined%
  \makeatother%
  \begin{picture}(1,0.71869135)%
    \lineheight{1}%
    \setlength\tabcolsep{0pt}%
    \put(0,0){\includegraphics[width=\unitlength,page=1]{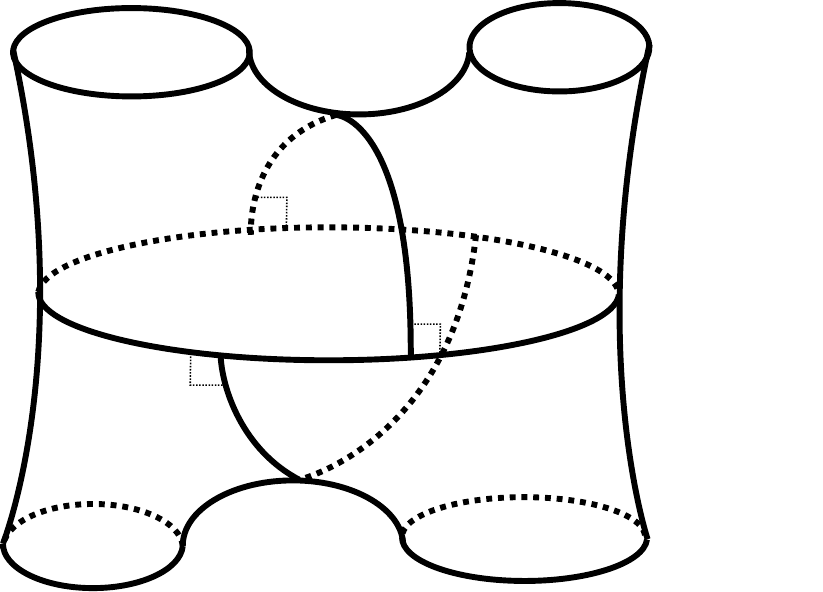}}%
    \put(0.48458769,0.49042554){\color[rgb]{0,0,0}\makebox(0,0)[lt]{\lineheight{1.25}\smash{\begin{tabular}[t]{l}$\eta_1$\end{tabular}}}}%
    \put(0.17701078,0.17316661){\color[rgb]{0,0,0}\makebox(0,0)[lt]{\lineheight{1.25}\smash{\begin{tabular}[t]{l}$\eta_2$\end{tabular}}}}%
    \put(0.59430568,0.23419472){\color[rgb]{0,0,0}\makebox(0,0)[lt]{\lineheight{1.25}\smash{\begin{tabular}[t]{l}$\gamma$\end{tabular}}}}%
    \put(0,0){\includegraphics[width=\unitlength,page=2]{fourholedlinked.pdf}}%
  \end{picture}%
\endgroup%
}}
\caption{Simple orthogeodesics $\eta_{1,2}$ to $\gamma$ with alternating endpoints in $\Sigma_{0,4}$.}
\label{fig:fourholedlinked}
}
\end{figure}
    
    \item The endpoints of $\eta_{1,2}$ do not alternate (see Figure~\ref{fig:fourholed}): Consider the disjoint arcs $\gamma_{1,2}$ in $\gamma$ with the same endpoints as $\eta_{1,2}$, respectively. Observe that $\ell(\eta_1) + \ell(\eta_2) \geq \ell(\gamma_1) + \ell(\gamma_2)$, as otherwise the concatenation of $\eta_1, \eta_2, \gamma\setminus \gamma_1, \gamma\setminus \gamma_2$ is a non-peripheral curve with smaller length than $\gamma$. Then for at least one index $i=1,2$ we must have $\ell(\eta_i)\geq \ell(\gamma_i)$. Without loss of generality we will assume that this is the case for $i=1$.
    
    Consider the region in the complement of $\eta_1$ that contains $\gamma_1$ and double it through $\eta_1$. This yields a pair of pants $P$ with a cuff $\alpha$ of length $2\ell(\gamma_1)$ and orthogeodesics $\eta_{1}$ (see Figure~\ref{fig:doubling}). As in the previous case we can apply the Collar Lemma~\cite[4.1.1.~Theorem]{Bus10:GeometryRiemannSurfaces} to the double of $P$ to obtain half collars around $\alpha$ of width $w$ satisfying $\sinh(w)\sinh(\ell(\gamma_1)) = 1$, which in turn implies $\sinh(w)\sinh(\eta_1)\geq 1$. But once again, as orthogeodesics are obstructions to half collars, we have that $\ell(\eta_1)\geq 2w$. The last two inequalities imply $\sinh(\frac12\ell(\eta_1))\sinh(\ell(\eta_1))\geq 1$, from where it follows that $\ell(\eta_1)\geq w_4\simeq 1.2188$, where $w_4$ is the positive constant satisfying $\sinh(w_4)\sinh(\frac12 w_4)=1$.
        \begin{figure}[H]
\centering{
\resizebox{100mm}{!}{\Huge{
\begingroup%
  \makeatletter%
  \providecommand\color[2][]{%
    \errmessage{(Inkscape) Color is used for the text in Inkscape, but the package 'color.sty' is not loaded}%
    \renewcommand\color[2][]{}%
  }%
  \providecommand\transparent[1]{%
    \errmessage{(Inkscape) Transparency is used (non-zero) for the text in Inkscape, but the package 'transparent.sty' is not loaded}%
    \renewcommand\transparent[1]{}%
  }%
  \providecommand\rotatebox[2]{#2}%
  \newcommand*\fsize{\dimexpr\f@size pt\relax}%
  \newcommand*\lineheight[1]{\fontsize{\fsize}{#1\fsize}\selectfont}%
  \ifx\svgwidth\undefined%
    \setlength{\unitlength}{394.75866603bp}%
    \ifx\svgscale\undefined%
      \relax%
    \else%
      \setlength{\unitlength}{\unitlength * \real{\svgscale}}%
    \fi%
  \else%
    \setlength{\unitlength}{\svgwidth}%
  \fi%
  \global\let\svgwidth\undefined%
  \global\let\svgscale\undefined%
  \makeatother%
  \begin{picture}(1,0.71869135)%
    \lineheight{1}%
    \setlength\tabcolsep{0pt}%
    \put(0,0){\includegraphics[width=\unitlength,page=1]{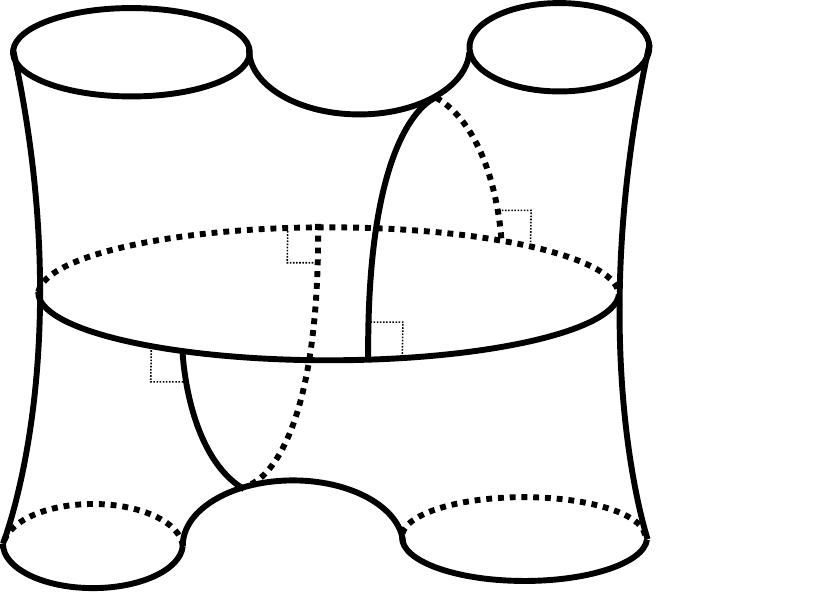}}%
    \put(0.49160585,0.46580544){\color[rgb]{0,0,0}\makebox(0,0)[lt]{\lineheight{1.25}\smash{\begin{tabular}[t]{l}$\eta_1$\end{tabular}}}}%
    \put(0.24418222,0.19619683){\color[rgb]{0,0,0}\makebox(0,0)[lt]{\lineheight{1.25}\smash{\begin{tabular}[t]{l}$\eta_2$\end{tabular}}}}%
    \put(0.59430568,0.23419472){\color[rgb]{0,0,0}\makebox(0,0)[lt]{\lineheight{1.25}\smash{\begin{tabular}[t]{l}$\gamma$\end{tabular}}}}%
  \end{picture}%
\endgroup%
}}
\caption{Simple orthogeodesics $\eta_{1,2}$ to $\gamma$ with non-alternating endpoints in $\Sigma_{0,4}$.}
\label{fig:fourholed}
}
\end{figure}
\end{itemize}

\end{proof}

We now reestate Theorem~\ref{thm:mainthm} for convenience, as Theorem~\ref{thm:thickannuli_bis}. 
\begin{theorem}
There exist universal constants
\[
w_0\simeq 0.40352, \qquad
c = \tfrac{2}{\pi-2\theta} \simeq \tfrac{2}{\pi-2\cdot 1.1778} \simeq \tfrac{2}{\pi-2.3556}, \qquad
\epsilon=f^{-1}(w_0)\approx 3.2282
\]
with the following property.  
Let $\Sigma$ be a closed hyperbolic surface, and let
$\mathcal{C}=\{\gamma_i\}_{1\leq i \leq n}$, for $n<3g-3$, be an $\epsilon$-short partial pants decomposition. Let $\gamma$ be
the shortest essential curve in a component of
$\Sigma \setminus \mathcal{C}$ that contains at least one simple essential
curve. Then $\gamma$ admits an embedded annular half collar
$A_{w_0}(\gamma)$ of width $w_0$. In particular, its extremal length satisfies
\[
  \EL(\gamma;\Sigma) \leq c\,\ell(\gamma;\Sigma).
\]
\label{thm:thickannuli_bis}
\end{theorem}

\begin{proof}[Proof of Theorem~\ref{thm:thickannuli_bis}]
We will assume that the component $\Sigma_0$ of $\Sigma\setminus\mathcal{C}$ containing $\gamma$ is not a three-holed sphere (if it is, there are no essential simple curves, and we are done). \\
Proceeding by contradiction, assume that there exists an orthogeodesic $\eta$ to $\gamma$ with $\ell(\eta)<2w_0$. If $\eta$ intersects $\mathcal{C}$ then we have that the distance between $\gamma$ and $\mathcal{C}$ is less than $w_0$. But since from Lemma \ref{lem:noshortorthogeodesics} we have that such distance is bounded below by $f(\epsilon)=w_0$ we have a contradiction. Hence we have that $\eta$ is disjoint from $\mathcal{C}$, which in other words means that $\eta$ is contained in $\Sigma_0$. By Lemma \ref{lem:ortho} then the only possibility is that $\gamma$ is separating in $\Sigma_0$ and the component of $\Sigma_0\setminus \gamma$ containing $\eta$ is a pair of pants. But again by Lemma \ref{lem:ortho} there is at least one component of $\Sigma_0\setminus\gamma$ where $\eta$ cannot exist, so we can take an embedded half collar $A_{w_0}(\gamma)$ in that component, finishing the contradiction.
By Theorem~\ref{thm:upperEL}, a half-cylinder of hyperbolic width $w_0$ around a geodesic $\gamma$ implies that $\EL(\gamma)\leq  \frac{\ell(\gamma)}{\pi/2-\theta}$, where $\sin\theta\cosh {w_0}=1$ 
for $\theta \approx 1.1778$, which concludes the statement.
\end{proof}

A priori, it's possible to have very short orthogeodesics going from a shortest curve to boundary parallel curves of a surface with boundary $\Sigma$.  However, taking the double of the surface we can still get bounds for the extremal length of short curves in the doubled surface, which induce bounds for the extremal length on $\Sigma$, as the following result, stated previously as Proposition~\ref{prop:boundary_systole}, shows. This provides a generalization of Theorem~\ref{thm:mainthm} to surfaces with boundary.  We reestate it here for convenience before proving it, and note that this result will not be used further.

\begin{figure}[H]
\centering{
\resizebox{80mm}{!}{\Huge{
\begingroup%
  \makeatletter%
  \providecommand\color[2][]{%
    \errmessage{(Inkscape) Color is used for the text in Inkscape, but the package 'color.sty' is not loaded}%
    \renewcommand\color[2][]{}%
  }%
  \providecommand\transparent[1]{%
    \errmessage{(Inkscape) Transparency is used (non-zero) for the text in Inkscape, but the package 'transparent.sty' is not loaded}%
    \renewcommand\transparent[1]{}%
  }%
  \providecommand\rotatebox[2]{#2}%
  \newcommand*\fsize{\dimexpr\f@size pt\relax}%
  \newcommand*\lineheight[1]{\fontsize{\fsize}{#1\fsize}\selectfont}%
  \ifx\svgwidth\undefined%
    \setlength{\unitlength}{369.73306491bp}%
    \ifx\svgscale\undefined%
      \relax%
    \else%
      \setlength{\unitlength}{\unitlength * \real{\svgscale}}%
    \fi%
  \else%
    \setlength{\unitlength}{\svgwidth}%
  \fi%
  \global\let\svgwidth\undefined%
  \global\let\svgscale\undefined%
  \makeatother%
  \begin{picture}(1,0.72644572)%
    \lineheight{1}%
    \setlength\tabcolsep{0pt}%
    \put(0,0){\includegraphics[width=\unitlength,page=1]{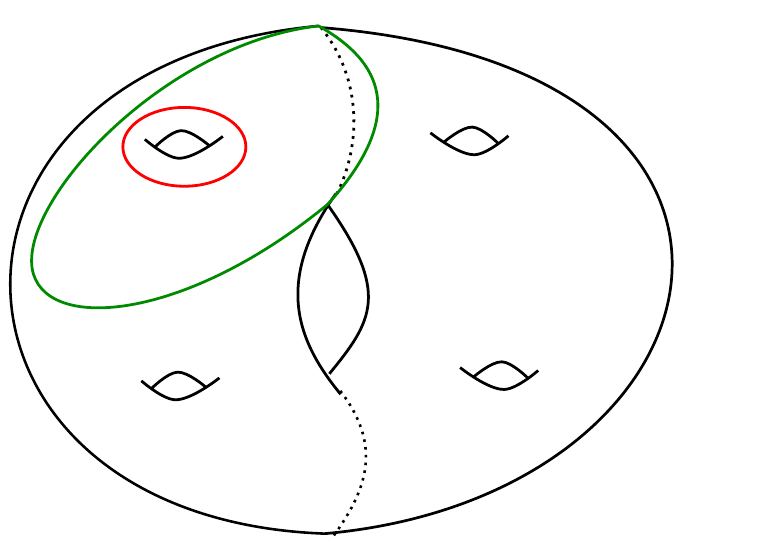}}%
    \put(0.53988038,0.48814182){\color[rgb]{0,0.54117647,0}\makebox(0,0)[t]{\lineheight{1.25}\smash{\begin{tabular}[t]{c}$\gamma$\end{tabular}}}}%
    \put(0.29086337,0.28128841){\color[rgb]{1,0,0}\makebox(0,0)[t]{\lineheight{1.25}\smash{\begin{tabular}[t]{c}$\tilde{\gamma}$\end{tabular}}}}%
    \put(0.05311122,0.62033618){\color[rgb]{0,0,0}\makebox(0,0)[lt]{\lineheight{1.25}\smash{\begin{tabular}[t]{l}$\Sigma$\end{tabular}}}}%
    \put(0.70038312,0.6551103){\color[rgb]{0,0,0}\makebox(0,0)[lt]{\lineheight{1.25}\smash{\begin{tabular}[t]{l}$D\Sigma$\end{tabular}}}}%
    \put(0,0){\includegraphics[width=\unitlength,page=2]{statement2.pdf}}%
  \end{picture}%
\endgroup%
}}
\caption{Relation between extremal length of a curve $\gamma$ in the surface with totally geodesic boundary $\Sigma$ and that of the curve in the double of the surface $D\Sigma$}
\label{fig:stat2}
}
\end{figure}

\begin{proposition}
Let $\Sigma$ be a hyperbolic surface with totally geodesic boundary.
Assume that each component of $\partial\Sigma$ has a collar of width $w$, and all such collars are disjoint. Let $\gamma$ be one of the shortest essential hyperbolic geodesics in $\Sigma$. Then 

\begin{equation}
\EL(\gamma;\Sigma)  \leq \frac{4\ell(\gamma)}{\pi-2\theta}\
\end{equation}
where $\theta\in[0,\pi/2)$ satisfies $\sin(\theta).\cosh(\min\lbrace w_0,w\rbrace)=1$ with $w_0\simeq 0.40352$.

\label{prop:boundary}
\end{proposition}
\begin{proof}
Let $\gamma$ be a shortest hyperbolic geodesic in $\Sigma$. Consider its homotopy class in the double of the surface $D\Sigma$. There are curves in the homotopy class of $\gamma$ in $D\Sigma$ that intersect $\sigma(\Sigma)$, i.e., the doubled copy of $\Sigma$. See for example the green curve in Figure~\ref{fig:stat2}, which we still denote by $\gamma$.
From such a representative, we can produce a new curve $\tilde{\gamma}$ also in the homotopy class of $\gamma$ in $D\Sigma$ as follows
\[
\tilde{\gamma}= \left( \gamma \cap \Sigma \right) \cup \left( \sigma^{*}(\gamma \cap \sigma(\Sigma)) \right). \]

Then, we can upper bound the extremal length of $\gamma$ in $\Sigma$ by the extremal length of the geodesic $\gamma$ in $D\Sigma$, as follows.
\begin{align*}
\EL(\gamma;D\Sigma)&=\sup_{\rho \in D\Sigma} \inf_{\gamma' \sim \gamma} \frac{\len(\gamma;\rho\Sigma)^2}{\Area(\rho)} \\
&\geq \sup_{\rho \in D\Sigma, \sigma^{*}(\rho)=\rho} \inf_{\gamma' \sim \gamma} \frac{\len(\gamma;\rho\Sigma)^2}{\Area(\rho)} \\
&=\sup_{\rho \in D\Sigma, \sigma^{*}(\rho)=\rho} \inf_{\gamma' \sim \gamma} \frac{\len(\tilde{\gamma};\rho\Sigma)^2}{\Area(\rho)} \\
&=\frac{1}{2}\EL(\gamma;\Sigma)
\end{align*}
where the first inequality is because we are restricting the metrics in the conformal class of $D\Sigma$, the next equality is because $\gamma$ and $\tilde{\gamma}$ are in the same homotopy class in $D\Sigma$, and the last equality is because any metric on $D\Sigma$ induces a metric on $\Sigma$ by restriction (which divides the area by a factor of two).

Note $\gamma$ has a half collar in $D\Sigma$  of width $w'=\min \{ w_0, w\}$ in $D\Sigma$. Indeed, suppose there is no such collar. Then, either
there exists an
orthogeodesic $\eta$ of length strictly smaller than $2w'$ reaching opposite sides of $\gamma$,  or we have orthogeodesics $\eta^+$ and $\eta^-$, reaching one side of $\gamma$ and both of length strictly smaller than $2w'$, reaching each side of $\gamma$.  In both cases, the orthogeodesics are contained in $\Sigma \subset D\Sigma$. But by Lemma~\ref{lem:ortho}, this is not possible, so we reach a contradiction. The conclusion then follows by Theorem~\ref{thm:upperEL}.
\end{proof}

\begin{remark}
Theorem~\ref{thm:thickannuli_bis} and Proposition~\ref{prop:boundary} also apply to surfaces with cusps.
Indeed, Lemma~\ref{lem:noshortorthogeodesics} and Theorem~\ref{thm:upperEL} do not use any assumptions on the existence of cusps and the proof in
Lemma~\ref{lem:ortho} follows in the case of cusps by allowing the auxiliary geodesics $\tilde{\gamma}$ to have length zero.
\end{remark}

        \begin{figure}[H]
\centering{
\resizebox{80mm}{!}{\Huge{
\begingroup%
  \makeatletter%
  \providecommand\color[2][]{%
    \errmessage{(Inkscape) Color is used for the text in Inkscape, but the package 'color.sty' is not loaded}%
    \renewcommand\color[2][]{}%
  }%
  \providecommand\transparent[1]{%
    \errmessage{(Inkscape) Transparency is used (non-zero) for the text in Inkscape, but the package 'transparent.sty' is not loaded}%
    \renewcommand\transparent[1]{}%
  }%
  \providecommand\rotatebox[2]{#2}%
  \newcommand*\fsize{\dimexpr\f@size pt\relax}%
  \newcommand*\lineheight[1]{\fontsize{\fsize}{#1\fsize}\selectfont}%
  \ifx\svgwidth\undefined%
    \setlength{\unitlength}{396.25427727bp}%
    \ifx\svgscale\undefined%
      \relax%
    \else%
      \setlength{\unitlength}{\unitlength * \real{\svgscale}}%
    \fi%
  \else%
    \setlength{\unitlength}{\svgwidth}%
  \fi%
  \global\let\svgwidth\undefined%
  \global\let\svgscale\undefined%
  \makeatother%
  \begin{picture}(1,0.33311424)%
    \lineheight{1}%
    \setlength\tabcolsep{0pt}%
    \put(0,0){\includegraphics[width=\unitlength,page=1]{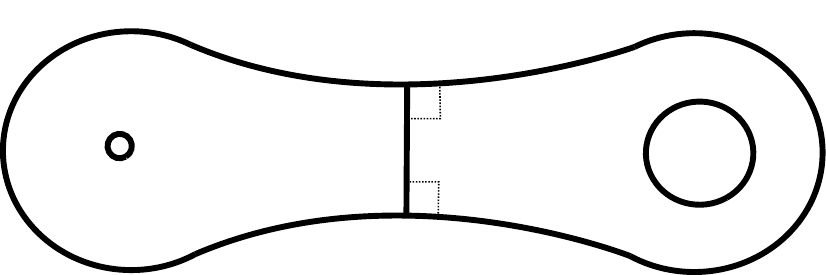}}%
    \put(0.79345245,0.31010581){\color[rgb]{0,0,0}\makebox(0,0)[lt]{\lineheight{1.25}\smash{\begin{tabular}[t]{l}$\gamma$\end{tabular}}}}%
    \put(0.49923554,0.1409309){\color[rgb]{0,0,0}\makebox(0,0)[lt]{\lineheight{1.25}\smash{\begin{tabular}[t]{l}$\eta$\end{tabular}}}}%
  \end{picture}%
\endgroup%
}}
\caption{Pair of pants with one geodesic boundary length going to 0.}
\label{fig:pantscusp}
}
\end{figure}


\section{Applications}

\subsection{Upper bounds for renormalized volume}
\label{subsec:RVol}
In order to give an answer to Question \ref{que:maldacena} we state Theorem 2.1 of \cite{VP19:UpperBounds}, which bounds $\VR$ of a handlebody in terms of extremal length of compressible curves. Recall that we say that a curve is \textit{compressible} if it bounds a topological disk in the 3-manifold. 
The following result is~\cite[Theorem~2.1]{VP19:UpperBounds}.

\begin{theorem}\label{thm:VREL}
Let $\Sigma$ be a closed Riemann surface of genus $g$, and let $\Gamma \coloneqq \lbrace \gamma_i\rbrace_{1\leq i\leq g-1}$ be a set of $g-1$ mutually disjoint, non-homotopic, simple curves of $\Sigma$ so that each component of $\Sigma-\Gamma$ has genus 1.  Let the sum of square roots of extremal lengths be denoted by
\[
L(\Sigma,\Gamma) \coloneqq \sum_{i=1}^{g-1} \sqrt{\EL(\gamma_i)}.
\]
Then for any Schottky manifold $M$ with boundary at infinity conformal to $\Sigma$ and where the curves of $\Gamma$ are compressible, we have that
\begin{equation}
    \VR(M) \leq L(\Sigma,\Gamma)^2 + \pi(g-1).
\end{equation}
Moreover, if we further assume that $L(\Sigma,\Gamma) \leq \sqrt{\pi(g-1)}$ then we have that

\begin{equation}
    \VR(M) \leq \pi(g-1)\left(3-\frac{\pi(g-1)}{L(\Sigma,\Gamma)^2}\right)
\end{equation}
which answers positively Maldacena's question if $L(\Sigma,\Gamma)\leq \sqrt{\frac{\pi(g-1)}{3}}$.
\end{theorem}
Applying our main Theorem~\ref{thm:mainthm} to \cite[Theorem 2.1]{VP19:UpperBounds} we obtain the following positive answer to Maldacena's question in terms of hyperbolic lengths of compressible curves.
\begin{corollary}\label{cor:vrsystole}
Let $w_0,\epsilon,c$ be the constants of Theorem~\ref{thm:mainthm}. Let $\Sigma$ be a closed Riemann surface and $\mathcal{C} = \lbrace \gamma_k\rbrace_{1\leq k \leq n}$ a collection of disjoint simple curves of length less than $\epsilon$ so that each $\gamma_k$ has an embedded half collar $A_{w_0}(\gamma_k)$ in $\Sigma$, so that the collection is suitable. Then we can extend $\mathcal{C}$ to a collection $\mathcal{C}'$ of disjoint curves so that each component of $\Sigma\setminus\mathcal{C}'$ has genus $1$, and so that for any Schottky manifold $M^3$ where all curves in $\mathcal{C}'$ are compressible we have the following upper bound for $\VR(M)$

\begin{equation}
    \VR(M) \leq \pi(g-1) + c\left(\sum_{i=1}^{g-1} \sqrt{\ell(\gamma_{k_i}, \Sigma)} \right)^2.
\end{equation}

Moreover, if $\sqrt{c} \left(\sum_{i=1}^{g-1} \sqrt{\ell(\gamma_{k_i}, \Sigma)} \right) \leq \sqrt{\pi(g-1)}$ then


\begin{equation}
    \VR(M) \leq \pi(g-1)\left(3-\frac{\pi(g-1)}{c\left(\sum_{i=1}^{g-1} \sqrt{\ell(\gamma_{k_i}, \Sigma)} \right)^2}\right).
\end{equation}
In particular $\VR\leq 0$ if $\sqrt{c}\sum_{i=1}^{g-1} \sqrt{\ell(\gamma_{k_i}, \Sigma)}\leq \sqrt{\frac{\pi(g-1)}{3}}$.
\label{cor:renvolume}
\end{corollary}
\begin{proof}
For each component of $\Sigma\setminus\mathcal{C}$ of genus $2$ we can add its shortest simple curves to $\mathcal{C}$ and define the collection as $\mathcal{C}'$. Since $\mathcal{C}$ is suitable and, in particular, each genus 2 component has only one boundary component, it follows that each component of $\Sigma\setminus\mathcal{C}'$ has genus $1$, and for each curve $\gamma$ in $\mathcal{C}'$ we have that $EL(\gamma) \leq c\ell(\gamma)$, by Theorem~\ref{thm:mainthm}. Then the bounds on $\VR(M)$ for a Schottky manifold $M$ where $\mathcal{C}'$ are compressible follow from Theorem \ref{thm:VREL}.
\end{proof}
\appendix

\section{Hyperbolic geometry}
\label{subsec:hypgeo}

We write some hyperbolic geometry formulas we will use in the sequel: the Collar Lemma \cite[Corollary~4.1.2]{Bus10:GeometryRiemannSurfaces}, Right-angled pentagon formula \cite[Theorem~2.3.4]{Bus10:GeometryRiemannSurfaces}, and Cosine formula \cite[Theorem~2.4.4]{Bus10:GeometryRiemannSurfaces}.

\begin{theorem}[Collar Lemma]
Let $\gamma, \delta$ be closed geodesics on $S$ which intersect transversally. Assume $\gamma$ is simple. Then

\[
\sinh(\frac12\ell(\gamma))\cdot \sinh(\frac12\ell(\delta))>1
\]
\label{eq:collar}
\end{theorem}
\begin{theorem}[Right-angled pentagon]
For any right-angled pentagon in the hyperbolic plane with consecutive sides $a,b,\alpha,c,\beta$ we have
\[
\cosh c = \sinh a \cdot  \sinh b
\]
By continuity, this formula still applies when  $c$ goes to $0$ and the picture degenerates into a trirectangle with acute angle $\phi=0$, satisfying

\[
1 = \sinh a \cdot  \sinh b
\]
\label{eq:rightanglepent}
\end{theorem}

\begin{theorem}[Cosine formula]
For any geodesic triangle in the hyperbolic plane with consecutive sides $a,b,c$ and angle $\gamma$ sustained between $a$ and $b$, we have

\[
\cosh c = \cosh a \cdot  \cosh b - \sinh a \cdot  \sinh b \cdot  \cos \gamma
\]
\label{eq:coslaw}
\end{theorem}

Finally, we state the upper bound by Maskit in \cite{Maskit85:ComparisonLengths} that can also be found in \cite[Equation~(3.5)]{BS94:PeriodMatrix}.

\begin{theorem}[Maskit's upper bound]
Given a simple curve $\gamma$ with a half-annular collar $A_w$ of hyperbolic width $w$ inside the surface $\Sigma$, we have
\[
  \EL(\gamma;\Sigma) \leq \frac{\ell(\gamma)}{\pi/2-\theta}
\]
where $0 < \theta < \frac{\pi}{2}$ is determined by

\[
  \cosh w \cdot \sin \theta=1.
\]
\label{thm:upperEL}
\end{theorem}

\bibliographystyle{hamsalpha}
\bibliography{main}

@book{MatsuzakiTaniguchi1998:HyperbolicManifolds,
  author    = {Katsuhiko Matsuzaki and Masahiko Taniguchi},
  title     = {Hyperbolic Manifolds and Kleinian Groups},
  series    = {Oxford Mathematical Monographs},
  publisher = {Clarendon Press},
  address   = {Oxford},
  year      = {1998},
  isbn      = {9780198500629},
  doi       = {10.1093/oso/9780198500629.001.0001}
}

@book{Marden2016:HyperbolicManifolds,
  author    = {Albert Marden},
  title     = {Hyperbolic Manifolds: An Introduction in 2 and 3 Dimensions},
  publisher = {Cambridge University Press},
  address   = {Cambridge},
  year      = {2016},
  isbn      = {9781107116740},
  doi       = {10.1017/CBO9781316337776}
}

@book {FB12:Primer,
    AUTHOR = {Farb, Benson and Margalit, Dan},
     TITLE = {A primer on mapping class groups},
    SERIES = {Princeton Mathematical Series},
    VOLUME = {49},
 PUBLISHER = {Princeton University Press, Princeton, NJ},
      YEAR = {2012},
     PAGES = {xiv+472},
      ISBN = {978-0-691-14794-9},
   MRCLASS = {57M50 (20F36 20F65 57M07 57N05)},
  MRNUMBER = {2850125},
MRREVIEWER = {Stephen\ P.\ Humphries},
}

@article{Krasnov2000,
  author        = {Krasnov, Kirill},
  title         = {Holography and Riemann Surfaces},
  journal       = {Advances in Theoretical and Mathematical Physics},
  volume         = {4},
  number         = {4},
  year           = {2000},
  pages          = {929--979},
  doi            = {10.4310/ATMP.2000.v4.n4.a5},
  eprint         = {hep-th/0005106},
  archivePrefix  = {arXiv},
  primaryClass   = {hep-th}
}

@article{SchlenkerWitten2022,
  author        = {Schlenker, Jean-Marc and Witten, Edward},
  title         = {No Ensemble Averaging Below the Black Hole Threshold},
  journal       = {Journal of High Energy Physics},
  volume        = {2022},
  number        = {7},
  pages         = {143},
  year          = {2022},
  doi           = {10.1007/JHEP07(2022)143},
  eprint        = {2202.01372},
  archivePrefix = {arXiv},
  primaryClass  = {hep-th}
}

@book{Bus10:GeometryRiemannSurfaces, title={Geometry and Spectra of Compact Riemann Surfaces}, DOI={10.1007/978-0-8176-4992-0}, author={Buser, Peter}, year={2010},  publisher = {Birkh{\"a}user/Springer}}

@book{Ahl10:ConformalInvariants, title={Conformal Invariants}, DOI={10.1090/chel/371}, author={Ahlfors, Lars}, publisher={AMS Chelsea Publishing}, place={Providence, Rhode Island},year={2010}}

@article{Maskit85:ComparisonLengths, title={Comparison of hyperbolic and extremal lengths}, volume={10}, DOI={10.5186/aasfm.1985.1042}, journal={Annales Academiae Scientiarum Fennicae. Series A. I. Mathematica}, author={Maskit, Bernard}, year={1985}, pages={381–386}}

@misc{MGT20:FromCurvesToCurrents,
      title={From curves to currents}, 
      author={D\'{i}dac Mart\'{i}nez-Granado and Dylan P. Thurston},
      year={2020},
      eprint={2004.01550},
      archivePrefix={arXiv},
      primaryClass={math.GT}
}

@article{FBMGVP:Bolza,
     AUTHOR = {Fortier Bourque, Maxime and Mart\'inez-Granado, D\'idac and
              Vargas Pallete, Franco},
     TITLE = {The extremal length systole of the {B}olza surface},
   JOURNAL = {Ann. H. Lebesgue},
  FJOURNAL = {Annales Henri Lebesgue},
    VOLUME = {7},
      YEAR = {2024},
     PAGES = {1409--1455},
      ISSN = {2644-9463},
   MRCLASS = {30F60 (31A15 57K20)},
  MRNUMBER = {4883873},
}

@article {LR11:Quasiconvex,
    AUTHOR = {Anna Lenzhen and Kasra Rafi},
     TITLE = {Length of a curve is quasi-convex along {T}eichm\"{u}ller space},
   JOURNAL = {Journal of Differential Geometry},
    VOLUME = {88},
      YEAR = {2011},
     PAGES = {267--295},
}

@article {Ker90:Asymptotics,
    AUTHOR = {Kerckhoff, Steven P.},
     TITLE = {The asymptotic geometry of {T}eichm\"{u}ller space},
   JOURNAL = {Topology},
  FJOURNAL = {Topology. An International Journal of Mathematics},
    VOLUME = {19},
      YEAR = {1980},
    NUMBER = {1},
     PAGES = {23--41},
      ISSN = {0040-9383},
   MRCLASS = {32G15 (30F10 57R30)},
  MRNUMBER = {559474},
MRREVIEWER = {C. Earle},
       DOI = {10.1016/0040-9383(80)90029-4},
       URL = {https://doi.org/10.1016/0040-9383(80)90029-4},
}

@article{Str84:QuadraticDiff, title={Quadratic Differentials}, DOI={10.1007/978-3-662-02414-0_2}, journal={Quadratic Differentials}, author={Strebel, Kurt}, year={1984}, pages={16–26}}

@ARTICLE{KPT17:ConformalEmbeddings,
       author = {{Kahn}, Jeremy and {Pilgrim}, Kevin M. and {Thurston}, Dylan P.},
        title = "{Conformal surface embeddings and extremal length}",
      journal = {arXiv e-prints},
         year = "2015",
        month = "Jul",
          eid = {arXiv:1507.05294},
        pages = {arXiv:1507.05294},
archivePrefix = {arXiv},
       eprint = {1507.05294},
 primaryClass = {math.CV},
       adsurl = {https://ui.adsabs.harvard.edu/abs/2015arXiv150705294K}
}

@article{Jenkins57:OnExistenceOfExtremalLength,
  title={On the Existence of Certain General Extremal Metrics},
  volume={66},
  DOI={10.2307/1969901},
  number={3},
  journal={The Annals of Mathematics},
  author={Jenkins, James A.},
  year={1957},
  pages={440}}

@article{Gro83:FillingManifolds, title={Filling Riemannian manifolds}, volume={18}, DOI={10.4310/jdg/1214509283}, number={1}, journal={Journal of Differential Geometry}, author={Gromov, Mikhael}, year={1983}, pages={1–147}}

@article{Sch93:ShortestGeodesic, title={Riemann surfaces with shortest geodesic of maximal length}, volume={3}, DOI={10.1007/bf01896258}, number={6}, journal={Geometric and Functional Analysis}, author={Schmutz, P.}, year={1993}, pages={564–631}}

@article{VP19:UpperBounds,
      AUTHOR = {Vargas Pallete, Franco},
     TITLE = {Upper bounds on renormalized volume for {S}chottky groups},
   JOURNAL = {Math. Ann.},
  FJOURNAL = {Mathematische Annalen},
    VOLUME = {392},
      YEAR = {2025},
    NUMBER = {1},
     PAGES = {733--750},
      ISSN = {0025-5831,1432-1807},
   MRCLASS = {57K20 (30F35)},
  MRNUMBER = {4887772},
       DOI = {10.1007/s00208-024-03086-2},
       URL = {https://doi.org/10.1007/s00208-024-03086-2},
}

@article{Min96:ExtremalLength,
author = "Minsky, Yair N.",
doi = "10.1215/S0012-7094-96-08310-6",
fjournal = "Duke Mathematical Journal",
journal = "Duke Math. J.",
month = "05",
number = "2",
pages = "249--286",
publisher = "Duke University Press",
title = "Extremal length estimates and product regions in {T}eichmüller space",
url = "https://doi.org/10.1215/S0012-7094-96-08310-6",
volume = "83",
year = "1996"
}

@article{Rodin74:ExtremalLength,
  title={The method of extremal length},
  volume={80},
  DOI={10.1090/s0002-9904-1974-13517-2},
  number={4},
  journal={Bulletin of the American Mathematical Society},
  author={Rodin, Burton},
  year={1974},
  pages={587-–607}}

@article{BS94:PeriodMatrix, title={On the period matrix of a {R}iemann surface of large genus (with an Appendix by J.H. Conway and N.J.A. Sloane)}, volume={117}, DOI={10.1007/bf01232233}, number={1}, journal={Inventiones Mathematicae}, author={Buser, P. and Sarnak, P.}, year={1994}, pages={27–56}}

@article {KS08,
    AUTHOR = {Krasnov, Kirill and Schlenker, Jean-Marc},
     TITLE = {On the renormalized volume of hyperbolic 3-manifolds},
   JOURNAL = {Comm. Math. Phys.},
  FJOURNAL = {Communications in Mathematical Physics},
    VOLUME = {279},
      YEAR = {2008},
    NUMBER = {3},
     PAGES = {637--668},
      ISSN = {0010-3616},
   MRCLASS = {53C65 (32G81 53C80)},
  MRNUMBER = {2386723},
       DOI = {10.1007/s00220-008-0423-7},
       URL = {https://doi.org/10.1007/s00220-008-0423-7},
}

@article {Witten98,
    AUTHOR = {Witten, Edward},
     TITLE = {Anti de {S}itter space and holography},
   JOURNAL = {Adv. Theor. Math. Phys.},
  FJOURNAL = {Advances in Theoretical and Mathematical Physics},
    VOLUME = {2},
      YEAR = {1998},
    NUMBER = {2},
     PAGES = {253--291},
      ISSN = {1095-0761},
   MRCLASS = {81T30 (81T60 83E30)},
  MRNUMBER = {1633012},
MRREVIEWER = {Douglas J. Smith},
       DOI = {10.4310/ATMP.1998.v2.n2.a2},
       URL = {https://doi.org/10.4310/ATMP.1998.v2.n2.a2},
}

@unpublished{Epstein,
	author = {Epstein, Charles L.},
	title = {Envelopes of horospheres and Weingarten surfaces in hyperbolic 3-space},
	note = {pre-print, Princeton University},
	year = {1984},
}

@article {Schlenker13,
    AUTHOR = {Schlenker, Jean-Marc},
     TITLE = {The renormalized volume and the volume of the convex core of
              quasifuchsian manifolds},
   JOURNAL = {Math. Res. Lett.},
  FJOURNAL = {Mathematical Research Letters},
    VOLUME = {20},
      YEAR = {2013},
    NUMBER = {4},
     PAGES = {773--786},
      ISSN = {1073-2780},
   MRCLASS = {57N10 (30F60 32G15)},
  MRNUMBER = {3188032},
MRREVIEWER = {Joonhyung Kim},
       DOI = {10.4310/MRL.2013.v20.n4.a12},
       URL = {https://doi.org/10.4310/MRL.2013.v20.n4.a12},
}

@article {BC15,
    AUTHOR = {Bridgeman, Martin and Canary, Richard D.},
     TITLE = {Renormalized volume and the volume of the convex core},
   JOURNAL = {Ann. Inst. Fourier (Grenoble)},
  FJOURNAL = {Universit\'{e} de Grenoble. Annales de l'Institut Fourier},
    VOLUME = {67},
      YEAR = {2017},
    NUMBER = {5},
     PAGES = {2083--2098},
      ISSN = {0373-0956},
   MRCLASS = {57M50 (30F40 30F50 35F45)},
  MRNUMBER = {3732685},
MRREVIEWER = {Hongbin Sun},
       URL = {http://aif.cedram.org/item?id=AIF_2017__67_5_2083_0},
}

@article {BBB19,
    AUTHOR = {Bridgeman, Martin and Brock, Jeffrey and Bromberg, Kenneth},
     TITLE = {Schwarzian derivatives, projective structures, and the
              {W}eil--{P}etersson gradient flow for renormalized volume},
   JOURNAL = {Duke Math. J.},
  FJOURNAL = {Duke Mathematical Journal},
    VOLUME = {168},
      YEAR = {2019},
    NUMBER = {5},
     PAGES = {867--896},
      ISSN = {0012-7094},
   MRCLASS = {30F40 (37F30)},
  MRNUMBER = {3934591},
       DOI = {10.1215/00127094-2018-0061},
       URL = {https://doi.org/10.1215/00127094-2018-0061},
}

@misc {VP17,
	author = {Vargas Pallete, Franco},
	title = {Additive continuity of the renormalized volume under geometric limits},
	note = {\href{https://arxiv.org/abs/1708.04009}{arXiv:1708.04009 [math.DG]}},
}

@article{CremaschiGiovanniniSchlenker,
title = {Filling Riemann surfaces by hyperbolic Schottky manifolds of negative volume},
journal = {Journal of Geometry and Physics},
volume = {217},
pages = {105628},
year = {2025},
issn = {0393-0440},
doi = {https://doi.org/10.1016/j.geomphys.2025.105628},
url = {https://www.sciencedirect.com/science/article/pii/S0393044025002128},
author = {Tommaso Cremaschi and Viola Giovannini and Jean-Marc Schlenker},
}

\end{document}